\documentclass[11pt]{article}

\usepackage{hyperref}
\usepackage{amsmath}
\usepackage{amssymb}
\usepackage{amsthm}
\usepackage{latexsym}
\usepackage{color}
\usepackage{graphicx}
\usepackage{appendix}
\usepackage{enumerate}
\usepackage[T1]{fontenc}
\usepackage[english]{babel}
\usepackage[table]{xcolor}

\DeclareSymbolFont{calletters}{OMS}{cmsy}{m}{n}
\DeclareSymbolFontAlphabet{\mathcal}{calletters}

\def\be{\begin{eqnarray}}
\def\ee{\end{eqnarray}}

\def\b*{\begin{eqnarray*}}
\def\e*{\end{eqnarray*}}

\newtheorem{Theorem}{Theorem}[part]

\newtheorem{Proposition}{Proposition}[part]

\newtheorem{Assumption}{Assumption}[part]
\newtheorem{Lemma}{Lemma}[part]
\newtheorem{Corollary}{Corollary}[part]
\newtheorem{Remark}{Remark}[part]
\newtheorem{Example}{Example}[part]

\makeatletter \@addtoreset{equation}{section}

\@addtoreset{Definition}{section}

\@addtoreset{Theorem}{section}

\@addtoreset{Proposition}{section}

\@addtoreset{Property}{section}

\@addtoreset{Assumption}{section}

\@addtoreset{Corollary}{section}

\@addtoreset{Lemma}{section}

\@addtoreset{Remark}{section}

\@addtoreset{Example}{section}

\newcommand{\No}[1]{\left\|#1\right\|}     
\newcommand{\abs}[1]{\left|#1\right|}     

\def \D{{\bf{D}}}
\def \E{\mathbb{E}}

\def \P{\mathbb{P}}

\def \R{\mathbb{R}}

\def\Ac{{\cal A}}

\def\Ic{{\cal I}}

\def\Lc{{\cal L}}
\def\Mc{{\cal M}}

\def\Ut{{\widetilde U}}

\def\cbf{{\bf c}}
\def\ybf{{\bf y}}

\def\Tr#1{{\rm Tr}\left[#1\right]}

\def \Frac{\displaystyle\frac}

\def\no{\noindent}
\def\x{\times}

\def\={\;=\;}
\def\.{\;.}

\def\eps{\varepsilon}

\def\reff#1{{\rm(\ref{#1})}}

\def\1{{\bf 1}}
\def \ep{\hbox{ }\hfill{ ${\cal t}$~\hspace{-5.1mm}~${\cal u}$   } }

\def \proof{{\noindent \bf Proof. }}
\def\ep{\epsilon}
\def\cbf{{\mathbf{c}}}

\def\ve{v^\epsilon}
\def\ep{\epsilon}
\def\eps{\epsilon}
\def\ue{u^\epsilon}

\def\b*{\begin{eqnarray*}}
\def\e*{\end{eqnarray*}}

 \def\normeL2#1{\left\|{#1}\right\|_{L^2}}
 
 \title{Homogenization and asymptotics for small transaction costs: the multidimensional case}

 \author{Dylan {\sc Possama\"{i}} \footnote{CEREMADE, Universit\'e Paris Dauphine, possamai@ceremade.dauphine.fr.}
 \and H. Mete {\sc Soner}\footnote{ETH (Swiss Federal Institute of Technology),
Zurich, and Swiss Finance Institute,
hmsoner@ethz.ch. Research partly supported by the
European Research Council under the grant 228053-FiRM,
the Swiss Finance Institute
and by the ETH Foundation.}
      \and Nizar {\sc Touzi}\footnote{CMAP, Ecole Polytechnique Paris, nizar.touzi@polytechnique.edu.
      Research supported by the Chair {\it Financial Risks} of the {\it Risk Foundation} sponsored by Soci\'et\'e
             G\'en\'erale, and
             the Chair {\it Finance and Sustainable Development} sponsored by EDF and Calyon. }}
             
 \date{Numerical results by: Lo\"{i}c {\sc Richier} and Bertrand {\sc Rondepierre}\footnote{Ecole Polytechnique Paris, loic.richier@polytechnique.edu, bertrand.rondepierre@polytechnique.edu.}
\\
\vspace{1.5em}
\today}

\begin{document}

 \maketitle

 \begin{abstract}
In the context of the multi-dimensional infinite horizon optimal consumption
investment problem with small proportional transaction costs, 
we prove an asymptotic expansion. Similar to the one-dimensional derivation in our accompanying paper \cite{st}, the first order term is expressed in terms of a singular ergodic control problem.
Our arguments are based on the theory of viscosity solutions and the techniques of homogenization which leads to a system of corrector equations. In contrast with the one-dimensional case, no explicit solution of the first corrector equation is available and we also prove the existence of a corrector and its properties. Finally, we provide some numerical results which illustrate the structure of the first order optimal controls.

\vspace{10mm}


\vspace{5mm}


\vspace{5mm}


\end{abstract}
\newpage

\section{Introduction}

We continue
our asymptotic analysis of
problems with small transaction costs
using the approach developed
in \cite{st} for problems
with only one stock.  
In this paper, we consider
the case of multi stocks with proportional
transaction costs.
The problem of investment and consumption in 
such a market  was first studied 
by Magill \& Constantinides  \cite{mc1976} and later by Constantinides \cite{con;86}. 
There is a large literature including the 
classical papers of Davis \& Norman \cite{da;no;90},
Shreve \& Soner \cite{ss}
and  Dumas \& Luciano \cite{dl1991}. 
We refer to our earlier paper \cite{st}
and to the recent book of Kabanov \& Safarian \cite{ks2009}
for other references and for more information.

\vspace{0.35em}
This problem is an important example
of a singular stochastic control
problem.  It is well known that
the related partial differential equation
contains a
gradient constraint.  
As such it is an interesting free boundary
problem.   
The main focus of this paper is on the analysis of the small transaction costs asymptotics.
 It is clear that in the limit of zero transaction costs, we recover the classical 
 problem of Merton \cite{merton} and the main interest is on the derivation 
 of the corrections of this obvious limit. 

\vspace{0.35em}
The asymptotic problem
is a  challenging problem which attracted 
considerable attention in the existing literature. The first rigorous proof in this direction was obtained in the appendix of \cite{ss}.
Later several rigorous results \cite{b2011,gms2011,js2004,r2004} and formal 
asymptotic results \cite{am2004,go2010,wh;wi;95} have been obtained. 
The rigorous results have been restricted to one space dimensions with 
the exception of the recent manuscript by Bichuch and Shreve \cite{bs2011}. 
As well known,
utility indifference price in this 
market is an important
approach as perfect hedging
is very costly as shown in \cite{ssc}.
Davis, Panas and Zariphopoulou \cite{dpz}
was first to study this approach with an exponential
utility function and the formal asymptotics 
was later developed in \cite{wh;wi;95}.

\vspace{0.35em}
In this paper, we use the 
techniques developed in \cite{st}.  
As in that paper, the main technique is the viscosity approach of Evans to homogenization \cite{ev1,ev2}.
This powerful method combined with the relaxed limits of 
Barles \& Perthame \cite{bp} provides the necessary tools. 
As well known, this approach has the advantage of using only a simple $L^\infty$ bound. 
In addition to \cite{bp,ev1,ev2}, the rigorous proof
utilizes several other techniques from the theory of viscosity solutions developed in the papers \cite{bp,fs89,fso,lsst,pst1,sszj,son93} for asymptotic analysis.

\vspace{0.35em}
For the classical problem of homogenization, we refer to the reader
to  the classical papers of Papanicolau \& Varadhan
\cite{pv79}
and of Souganidis \cite{sou3}.
However, we emphasize that the problem studied here 
does not even look like a typical homogenization problem of the existing literature. Indeed, the dynamic programming equation, which is the starting point of our asymptotic analysis,
does not include oscillatory variables of the form $x/\eps$ or in general a fast varying ergodic quantity.
The fast variable appears after a change of variables that includes the difference with the optimal strategy at the effective level problem. Hence, the fast variable actually depends upon
the behavior of the limit problem, which is a  novel viewpoint for 
homogenization where usually the fast variables are built into the equation 
from the beginning and one does not have to construct them from
limits.  This ergodic variable, in general,  
does not have to be periodic.
In recent studies,
deep techniques combining ergodic theorems
and difficult parabolic estimates were used
to study these more general cases.  We refer the reader to
Lions \& Souganidis \cite{sou2} for the almost-periodic case
and Caffarelli \& Souganidis \cite{sou1} for the case in random media 
and to the references therein.

\vspace{0.35em}
As in our accompanying paper \cite{st}, the formal asymptotic analysis leads to a system of corrector equations related
to an ergodic optimal control problem similar to the 
monotone follower \cite{bsw}.
However, in contrast with the one-dimensional case studied in \cite{st}, no explicit solution is available for this singular ergodic control problem. This is the main difficulty that we face in the present multi-dimensional setting. 
We use the recent analysis of Hynd \cite{hyn,hyn2}
to analyze this multidimensional
problem and obtain a $C^{1,1}$
unique solution of the corresponding
``eigenvalue'' problem
satisfying a precise growth condition.
This characterization is sufficient
to carry out the asymptotic analysis.
On the other hand,
the regularity of the free boundary
is a difficult problem and we do not 
study it in this paper.  We refer
the reader to \cite{ss1,ss2}
for such analysis in a similar problem.

\vspace{0.35em}

The paper is organized as follows. Section \ref{sect:general setting} provides a quick review of the infinite horizon optimal consumption-investment problem under transaction costs, an recalls the formal asymptotics, as derived in \cite{st}. Section \ref{sect:main result} collects the main results of the paper. The next section is devoted to the numerical experiments which illustrate the nature of the first order optimal transfers. The rigorous proof of the first order expansion is reported in Section \ref{sect:convergence}. In particular, the proof requires some wellposedness results of the first corrector equation which are isolated in Section \ref{sect:first corrector}.

\vspace{3mm}

\no {\bf Notations:}\quad Throughout the paper, we denote by $\cdot$ the Euclidean scalar product in $\R^d$, and by $(e_1,\ldots,e_d)$ the canonical basis of $\R^d$. We shall be usually working on the space $\R\x\R^d$ with first component enumerated by $0$. The corresponding canonical basis is then denoted by $(e_0,\ldots,e_d)$. We denote by $\Mc_{d}(\R)$ the space of $d\times d$ matrices with real entries, and 
by $^{\rm{T}}$ the transposition of matrices. We denote by $B_r(x)$ the open ball of radius $r>0$ centered at $x$, and $\overline{B}_r(x)$ the corresponding closure.

\section{The general setting}
\label{sect:general setting}

In this section, we briefly review  the infinite horizon optimal consumption-investment problem under transaction costs and recall the formal asymptotics, derived in  \cite{st}.
These calculations are the starting point of our analysis.

\subsection{Optimal consumption and investment under proportional transaction costs}

The financial market consists of a non-risky asset $S^0$ and $d$ risky assets with price process $\{S_t=(S^1_t, \ldots,S^d_t),t\ge 0\}$ 
given by the stochastic differential equations (SDEs),
 $$
 \frac{dS^0_t}{S^0_t}
 =
 r(S_t)dt,~~
 \frac{dS^i_t}{S^i_t}
 =
 \mu^i(S_t) dt + \sum_{j=1}^d \sigma^{i,j}(S_t) dW^j_t,
 ~~1\le i\le d,
 $$
where $r:\R^d\to \R_+$ is the instantaneous interest rate and $\mu:\R^d \to \R^d$, $\sigma:\R^d\to\Mc_{d}(\R)$ are the coefficients of instantaneous mean return and volatility, satisfying the standing assumptions:
 \begin{equation*}
 r, \mu, \sigma \mbox{ are bounded and Lipschitz, and }
 (\sigma \sigma^{T})^{-1}~~\mbox{is bounded.}
 \end{equation*}
In particular, this guarantees the existence and
the  uniqueness of a strong solution to the above 
stochastic differential equations (SDEs).

\vspace{0.35em}
The portfolio of an investor is represented by the dollar value $X$ invested in the non-risky asset and the vector process $Y=(Y^1,\ldots,Y^d)$ of 
the value of the positions in each risky asset. These state variables are controlled by the choices of the total amount of transfers $L^{i,j}_t$, $0\le i,j\le d$, from the $i$-th to the $j$-th asset cumulated up to time $t$. Naturally, the control processes $\{L^{i,j}_t,t\ge 0\}$ are defined
as c\`ad-l\`ag, nondecreasing, adapted processes with $L_{0^-}=0$ and $L^{i,i}\equiv 0$. 

\vspace{0.35em}
In addition to the trading activity, the investor consumes at a rate determined by a nonnegative progressively measurable process $\{c_t,t\ge 0\}$. Here $c_t$ represents the rate of consumption in terms of the non-risky asset $S^0$. Such a pair $\nu:=(c,L)$ is called a {\em consumption-investment strategy}. For any initial position $(X_{0^-},Y_{0^-})=(x,y)\in\R\x\R^{d}$, the portfolio positions of the investor are given by the following state equation
 $$
 dX_t
 =
 \big(r(S_t)X_t-c_t\big)dt
 +\mathbf{R}^0(dL_t),
 ~~\mbox{and}~~
 dY^i_t
 =
  Y^i_t\;\frac{dS^i_t}{S^i_t}
 +\mathbf{R}^i(dL_t),
 ~~i=1,\ldots,d,
 $$
where 
$$
\mathbf{R}^i(\ell) 
:= 
\sum_{j=0}^{d}
 \big(\ell^{j,i}-(1+\eps^3\lambda^{i,j})\ell^{i,j}\big),
 ~~i=0,\ldots,d,
 \text{ for all } 
 \ell\in\Mc_{d+1}(\R_+),
$$
is the change of the investor's position in the $i-$th asset induced by a transfer policy $\ell$, given a structure of proportional transaction costs $\eps^3\lambda^{i,j}$ for any transfer from asset $i$ to asset $j$. Here, $\eps>0$ is a small parameter, $\lambda^{i,j}\ge 0$, $\lambda^{i,i}=0$, for all $i,j=0,\ldots,d$, and the scaling $\ep^3$ is chosen to state the expansion results simpler.

\vspace{0.35em}
Let $(X,Y)^{\nu,s,x,y}$ denote the controlled state process. A consumption-investment strategy $\nu$ is said to be {\em admissible} for the initial position $(s,x,y)$, if the induced state process satisfies the solvency condition $(X,Y)^{\nu,s,x,y}_t\in K_\eps,$ for all $t\ge 0$, $\P-$a.s., where the solvency region is defined by:
 \begin{align*}
 K_\eps
 :=
 \left\{(x,y)\in\R \times \R^{d}: \ (x,y)+\mathbf{R}(\ell) \in \R_+^{1+d} \text{ for some }\ell\in\Mc_{d+1}(\R_+) \right\}.
 \end{align*}
The set of admissible strategies is denoted by $\Theta^\ep(s,x,y)$. For given initial positions $S_{0}=s \in \R_+^d$, $X_{0^-}=x \in \R$, $Y_{0^-}=y \in \R^d$, the consumption-investment problem is the following maximization problem,
 \b*
 \ve(s,x,y)
 &:=&
 \sup_{(c,L) \in \Theta^\ep(s,x,y)}\
 \E\left[\int_0^\infty\ e^{-\beta t}\ U(c_t)dt \right],
 \e*
where 
$U:(0,\infty)\mapsto \R$ is a utility function. We assume that $U$ is $C^2$, increasing, strictly concave, and we denote its convex conjugate by,
 \b*
 \Ut(\tilde c)
 &:=& 
 \sup_{c>0} \big\{U(c)-c\tilde c\big\},
 \qquad
 \tilde c\in\R.
 \e*

\subsection{Dynamic programming equation}

The dynamic programming equation corresponding to the singular stochastic control problem $v^\eps$ involves the following differential operators. Let:
\be\label{Lc}
 \Lc
 &:=&
 \mu\cdot\left(\D_s +\D_y\right)
 +r\D_x
 +\frac12\mbox{Tr}\left[\sigma \sigma^{\rm T}\left(\D_{yy}+\D_{ss}+2 \D_{sy}\right) \right],
 \ee
and for $i,j=1,\ldots,d$,
 \b*
 &\D_x:=x\Frac{\partial}{\partial x},
 ~~\D_s^i:=s^i\Frac{\partial}{\partial s^i},
 ~~\D_y^i:=y^i\Frac{\partial}{\partial y^i},&
 \\
 &\D_{ss}^{i,j}:=s^is^j\Frac{\partial^2}{\partial s^i\partial s^j},
 ~~\D_{yy}^{i,j}:=y^iy^j\Frac{\partial^2}{\partial y^i\partial y^j},
 ~~\D_{sy}^{i,j}:=s^iy^j\Frac{\partial^2}{\partial s^i\partial y^j},&
 \e*
  $\D_s=(\D_{s}^i)_{1\le i\le d}$, $\D_y=(\D_y^i)_{1\le i\le d}$, 
  $\D_{yy}:=(\D_{yy}^{i,j})_{1\le i,j\le d}$, $\D_{ss}:=(\D_{ss}^{i,j})_{1\le i,j\le d}$, 
  $\D_{sy}:=(\D_{sy}^{i,j})_{1\le i,j\le d}$.
Moreover, for a smooth scalar function $(s,x,y)\in \R^d_+\times\R\times\R^d\longmapsto \varphi(x,y)$, we set
 \b*
 \varphi_{x}:=\frac{\partial\varphi}{\partial x}\;\in\R,
\qquad
 \varphi_{y}:=\frac{\partial\varphi}{\partial y}\;\in\R^d.
 \e*

\begin{Theorem}
\label{t.dpp}
Assume that the value function $\ve$ is locally bounded. 
Then, $\ve$ is a viscosity solution of the dynamic programming equation in 
$\R_+^d\x K_\ep$,
 \begin{equation} \label{e.dpp}
 \min_{0\le i,j\le d}
 \left\{\ \beta \ve - \Lc \ve - \Ut(\ve_{x})\ , 
        \ \Lambda_{i,j}^\eps\cdot(\ve_x,\ve_y) \ 
 \right\} 
 =
 0, \ \Lambda_{i,j}^\eps
 :=
 e_i-e_j+\eps^3\lambda^{i,j}\;e_i,
 \ 0\le i,j\le d.
 \end{equation}
Moreover, $\ve$ is concave in $(x,y)$ and converges to the Merton value function 
$v:=v^0$, as $\eps>0$ tends to zero.
\end{Theorem}

The Merton value function $v=v^0$ corresponds to the limiting case $\eps=0$ where the transfers between assets are not subject to transaction costs. Our subsequent analysis assumes that $v$ is smooth, which can be verified under slight conditions on the coefficients. In this context, we recall that $v$ can be characterized as the unique classical solution (within a convenient class of functions) of the corresponding HJB equation:
$$
 \beta v - rz v_z - \Lc^0 v - \Ut(v_z)
 -\sup_{\theta\in\R^d}\Big\{\theta\cdot\big((\mu-r\1_d) v_z
                                             +\sigma \sigma^{\rm T}\D_{sz} v
                                       \big)
                              +\frac12|\sigma^{\rm T}\theta|^2 v_{zz}
                        \Big\}
 =
 0,
$$
where $\1_d:=(1,\ldots,1)\in \R^d$, $\D_{sz}:=\frac{\partial}{\partial z}\D_s$, and
 \be\label{Lc0}
 \Lc^0
 &:=&
 \mu\cdot\D_s
 +\frac12\mbox{Tr}\big[\sigma \sigma^{\rm{T}} \D_{ss}\big].
 \ee
The optimal consumption and positioning in the various assets are defined by the functions $\cbf(s,z)$ and $\ybf(s,z)$ 
obtained as the maximizers of the Hamiltonian:
 \begin{align*}
 \cbf(s,z)
 &:=
 -\Ut^{\prime}\left(v_z(s,z)\right)
 =\left(U^{\prime}\right)^{-1}\left(v_z(s,z)\right)
 ,\\
 -v_{zz}(s,z)\sigma\sigma^{\rm T}(s)\ybf(s,z)
 &:=
 (\mu-r\1_d)(s) v_z(s,z)
 +\sigma \sigma^{\rm T}(s)\D_{sz} v(s,z)
 \ \ \mbox{for}\ \ 
 s\in\R_+^d,\ \ z\ge 0.
 \end{align*}

\subsection{Formal Asymptotics}
 
 Here, we  recall the formal first order expansion derived in \cite{st}:
\begin{equation}\label{goal-expansion}
v^\eps(s,x,y)
=
v(s,z)-\eps^2u(s,z)+\circ(\eps^2),
\end{equation} 
where $u$ is solution of the second corrector equation:
\b*
\Ac u
:=\beta u 
 -\Lc^0 u
 - \big(rz + \ybf\cdot(\mu-r\1_d) - \cbf\big) 
   u_z
 -\frac12 |\sigma^{\rm T} \ybf|^2 \,u_{zz}
 -\sigma\sigma^{\rm T}\ybf\cdot\D_{sz} u
=
a,
\e*
and the function $a$ is the second component of the solution $(w,a)$ of the first corrector equation:
\b*
&&\underset{0\leq i,j\leq d}{\max}\max\left\{\frac{\abs{ \sigma(s)\xi}^2}{2}v_{zz}(s,z)
-\frac12\Tr{\alpha\alpha^T(s,z) w_{\xi\xi}(s,z,\xi)}+a(s,z)\ ;\right.\\
&&\left.\hspace{150pt}-\lambda^{i,j}v_z(s,z)+ \frac{\partial w}{\partial \xi_i}(s,z,\xi)
- \frac{\partial w}{\partial \xi_j}(s,z,\xi) \right\}=0,
\end{eqnarray*}
where $\xi\in\R^d$ is the dependent variable, while $(s,z) \in (0,\infty)^d \times \R^+$ are fixed, and the diffusion coefficient is given by
$$
\alpha(s,z):=\left[\left(I_d-\ybf_z(s,z)\1_d^T\right)
{\rm{diag}}[\ybf(s,z)]-\ybf_s^T(s,z){\rm{diag}}[s]\right]\sigma(s).
$$

The expansion \reff{goal-expansion} was proved rigorously in \cite{st} in the one-dimensional case, with a crucial use of the explicit solution of the first corrector equation in one space dimension. Our objective in this paper is to show that the above expansion is valid in the present $d-$dimensional framework where no explicit solution of the first corrector equation is available anymore. 

\vspace{0.35em}
We finally recall the from \cite{st} the following normalization. Set
$$
\eta(s,z):=-\frac{v_z(s,z)}{v_{zz}(s,z)},\ 
\rho:=\frac{\xi}{\eta(s,z)},\ \overline{w}(s,z,\rho)
:=\frac{w(s,z,\eta(s,z)\rho)}{\eta(s,z)v_z(s,z)},
$$ 
$$
\overline{a}(s,z):=\frac{a(s,z)}{\eta(s,z)v_z(s,z)}, \ 
\bar{\alpha}(s,z):=\frac{\alpha(s,z)}{\eta(s,z)},
$$
so that the corrector equations with variable 
$\rho \in \R^d$ have the form,
\begin{eqnarray}
\nonumber &&\underset{0\leq i,j\leq d}{\max}\max\left\{\frac{\abs{ \sigma(s)\rho}^2}{2}-
\frac12\Tr{\bar\alpha\bar\alpha^T(s,z) \overline{w}_{\rho\rho}(s,z,\rho)}+\overline a(s,z)\ ;\right.
\\
&&\left.\hspace{150pt}-\lambda^{i,j}+ \frac{\partial \overline{w}}{\partial \rho_i}(s,z,\rho)- 
\frac{\partial \overline{w}}{\partial \rho_j}(s,z,\rho) \right\}=0
\label{eq:corrector1}\\
&&\hspace{10pt}\mathcal Au(s,z)=v_z(s,z)\eta(s,z)\overline{a}(s,z).
\label{eq:corrector2}
\end{eqnarray}

Notice that we will always use the normalization $\overline{w}(s,z,0)=0$.

\section{The main results}
\label{sect:main result}

\subsection{The first corrector equation}

We first state the existence and uniqueness of a solution to the first corrector equation \reff{eq:corrector1}, within a convenient class of functions. We also state the main properties of the solution which will be used later. We will often make use of the language of ergodic control theory, and say 
that $\overline{w}$ is a solution of \reff{eq:corrector1} with eigenvalue $\overline{a}$. 

\vspace{0.35em}
Consider the following closed convex subset of $\mathbb R^d$, and the corresponding support function
 \b*
 C
 :=
 \left\{\rho\in\mathbb R^d:
        -\lambda^{j,i}\leq \rho_i-\rho_j\leq \lambda^{i,j}, 
        \ 0\leq i,j\leq d
 \right\},
 &\delta_C(\rho):=\underset{u\in C}{\sup}\ u \cdot \rho,& 
 \rho\in\mathbb R^d,
\e*
with the convention that $\rho_0=0$. Then, $C$ is a bounded convex polyhedral containing $0$, which corresponds to the 
intersection of $d(d+1)$ hyperplanes. Furthermore, the gradient constraints of the corrector equation \reff{eq:corrector1} is exactly that $D\overline{w}\in C$. In this respect, \reff{eq:corrector1} is closely related to the variational inequality studied by Menaldi et al. \cite{mrt} and Hynd \cite{hyn}, where the gradient is restricted to lie in the closed unit ball of $\mathbb R^d$ (we also refer to \cite{hyn2} for related variational inequalities with general gradient constraints).

\vspace{0.35em}
The wellposedness of the first corrector equation is stated in the following two results. Since the variables $(s,z)$ are frozen in the equation \reff{eq:corrector1}, we omit the dependence on them. 
 
\begin{Theorem}\label{prop.uni}{\rm (First corrector equation: comparison)}
Suppose $w_1$ is a viscosity subsolution of \reff{eq:corrector1} with eigenvalue $a_1$ and that $w_2$ is a viscosity supersolution of \reff{eq:corrector1} with eigenvalue $a_2$. Assume further that
$$
\underset{\abs{\rho}\rightarrow+\infty}{\overline\lim}\ \frac{w_1(\rho)}{\delta_C(\rho)} 
\leq1\leq \underset{\abs{\rho}\rightarrow+\infty}{\underline\lim}\ \frac{w_2(\rho)}{\delta_C(\rho)}.
$$
Then, $a_1\leq a_2$. 
\end{Theorem}

The proof of this result is given in Section \ref{sect:comparison w}.

\begin{Theorem}\label{prop}{\rm (First corrector equation: existence)}
There exists a solution $\overline{w} \in C^{1,1}$ of the 
equation \reff{eq:corrector1} with eigenvalue $\overline a$, satisfying the growth condition $\lim_{|\rho|\to\infty}(\overline{w}/\delta_C)(\rho)=1$. Moreover,
\\
$\bullet$ $\overline{w}$ is convex and positive.
\\
$\bullet$ The set $\mathcal O_0:=\left\{\rho\in\mathbb R^d,
	\ D\overline{w}(\rho)\in\text{$\rm{int}$}(C)\right\}$
	 is open and bounded, $\overline{w}\in C^\infty(\mathcal O_0)$ and
	\b*
	\overline{w}(\rho)=\underset{y\in\overline{\mathcal O}_0}{\inf}\left\{\overline{w}(y)+\delta_C(\rho-y)\right\},
	\ \mbox{for all}&
	\rho\in\R^d.
	\e*
$\bullet$	There is a constant $M>0$ such that $0\leq D^2\overline{w}(\rho)\leq M1_{\overline{\mathcal O}_0}(\rho)$ for a.e. $\rho\in\mathbb R^d$.
\end{Theorem}

The proof of this result will be reported in Section \ref{sect:existence w}. 

\begin{Remark}\label{rem.rep}
{\rm{
As already pointed out in \cite{st}, the corrector equation \reff{eq:corrector1} 
is related to the dynamic programming equation of an ergodic control problem \cite{bor,bor1}. More precisely, let $(M^{i,j}_t)_{t\geq 0}$ be non-decreasing control processes for $0\leq i,j\leq d$, and define the control process $\rho$ as the solution of the following SDEs
 $$
 \rho^i_t
 =
 \rho_0^i+\sum_{j=1}^d\bar{\alpha}^{i,j}B_t^j
 +\sum_{j=0}^d\left(M^{j,i}_t-M^{i,j}_t\right).
 $$
Then, the ergodic control problem is given by,
 $$
 \overline a:=\underset{M}{\inf}\ J(M),
 ~~
 J(M)
 :=
 \underset{T\rightarrow+\infty}{\overline{\lim}}\frac1T
 \mathbb E\Big[\frac12\int_0^T
               |\sigma\rho_t|^2dt+\sum_{i,j=0}^d\lambda^{i,j}M_T^{i,j}
          \Big].
 $$
However, it is not clear whether the corresponding potential function, which is the candidate for a solution to \reff{eq:corrector1}, satisfies the growth condition $\overline{w}\underline{+\infty}{\sim}\delta_C$. For this reason, we cannot use this probabilistic representation, and our analysis of the equation \reff{eq:corrector1} follows the PDE approach of Hynd \cite{hyn}.}}
\qed
\end{Remark}

In the one-dimensional context $d=1$, the unique solution of the first corrector equation is easily obtained in explicit form in \cite{st}. For higher dimension $d\ge 2$, in general, no such explicit expression is available anymore. The following example illustrates a particular structure of the parameters of the problem which allows to obtain an explicit solution.
\begin{Example}
\label{example1}
{\rm{Suppose
the equation is of the form,
\begin{equation}
\label{e.what}
\underset{0\leq i\leq d}
{\max}\max\left\{-\frac{c^*_1\abs{\rho}^2}{2}
-\frac{c^*_2}{2}\Delta \widehat{w}(\rho)+\widehat{a}\ ,\ -\widehat{\lambda}^{i}
+ \frac{\partial \widehat{w}}{\partial \rho_i}(\rho)\ ,\ -\tilde{\lambda}^{i}
- \frac{\partial \widehat{w}}{\partial \rho_i}(\rho)
\right\}=0,
\end{equation}
with  given positive constants
$c^*_1, c^*_2, \widehat{\lambda}^{i}, \tilde{\lambda}^{i}$
 and  the normalization $\widehat{w}(0)=0$.
This corresponds to the case when
$\sigma$ and $\bar{\alpha}$ are multiples of  the 
identity matrix and  we are only allowed
to make transactions to and from the cash account, i.e., when $\lambda^{i,j}=\infty$ as soon as $i$ and $j$ are both different from $0$. Then, the unique solution $\widehat{w}$ 
is given as,
$$
\widehat{w}(\rho) = \sum_{i=1}^{d}\ \widetilde{w}_i(\rho_i),
$$
where $\widetilde{w}_i$ is the explicit solution of 
the one dimensional problem constructed in \cite{st}.
Moreover, $\widehat{a}$ is an explicit constant
independent of $z$.

\vspace{0.35em}
Notice that for the corrector equations, $\sigma$ and $\bar{\alpha}$ cannot be specified independent of each other.
However, the above explicit solution will be used as an upper bound for the unique solution of the corrector equation. 
}}
\end{Example}

\subsection{Assumptions}

This section collects all assumptions which are needed for our main results, and comments on them.

\subsubsection{Assumptions on the Merton value function}

\begin{Assumption}[Smoothness]\label{assump2}
$v$ and $\theta$ are $C^2$ in $(0,\infty)^{d+1}$, $v_z>0$, $\ybf^i_z>0$, $i\le d$, on $(0,\infty)^{d+1}$, and there exist $c_0,c_1>0$ such that
	\b*
	c_0\leq \ybf_z\cdot \1_d\le 1-c_0
	&\mbox{and}&
	\alpha\alpha^{\rm T}\ge c_1 I_d,
	~~\mbox{on}~~(0,\infty)^{d+1}.
	\e* 
\end{Assumption}

Notice that this Assumption is verified in the case of the Black-Scholes model with power utility.

\subsubsection{Local boundedness}

As in \cite{st}, we define
\be\label{ubareps}
\bar{u}^\eps(s,x,y)
&:=&
\frac{v(s,z)-v^\eps(s,x,y)}{\eps^2},\ s\in\mathbb R_+^d, ~~(x,y)\in K_\eps.
\ee
Following the classical approach of Barles and Perthame, we introduce the relaxed semi-limits
$$
u^*(s,x,y)
:=
\underset{(\eps,s',x',y')\rightarrow (0,s,x,y)}{\overline{\lim}}\bar{u}^\eps(s',x',y'),\ u_*(s,x,y)
:=\underset{(\eps,s',x',y')\rightarrow (0,s,x,y)}{\underline{\lim}}\bar{u}^\eps(s',x',y').$$

The following assumption is verified in Lemma \ref{l.bound} below, in the power utility context with constant coefficients.

\begin{Assumption}[Local bound]\label{assump1}
The family of functions $\bar{u}^\eps$ is locally uniformly bounded from above.
\end{Assumption}

This assumption states that for any $(s_0,x_0,y_0) \in (0,\infty)^d\x\R\x\R^d$ with $x_0+y_0\cdot \1_d>0$, 
there exist $r_0=r_0(s_0,x_0,y_0)>0$ and $\ep_0=\ep_0(s_0,x_0,y_0)>0$ so that
 \begin{equation} 
 \label{e.bound}
 b(s_0,x_0,y_0):= 
 \sup\{\ \ue(s,x,y)\ :\ (s,x,y) \in B_{r_0}(s_0,x_0,y_0),\
 \ep \in (0,\ep_0]\ \} <\infty,
 \end{equation}
where $B_{r_0}(s_0,x_0,y_0)$ denotes the open ball with radius $r_0$, centered at $(s_0,x_0,y_0)$.

\vspace{0.35em}
The following result verifies the local boundedness under a natural hypothesis. We observe that this is just one possible set of assumptions, and the proof of the lemma can be modified to obtain the same result under other conditions.

\begin{Lemma}\label{l.bound}
Suppose that $U$ is a power utility, and $0\le \lambda^{0,j},\lambda^{j,0}<\infty$, $\lambda^{i,j}=\infty$ for all $1\le i\neq j\le d$. Then, Assumption \ref{assump1} holds.
\end{Lemma}
 
\vspace{0.35em}
\proof
The homotheticity of the utility function
implies that both the Merton value function
$v$ and the solution $u$ of the
second corrector equation are
homothetic as well.
For a large positive constant $K$ to be chosen later,
set
$$
V^{\eps,K}(z,\xi) := v(z) - \eps^2 K u(z) - \eps^4 W(z,\xi),
$$
where
$W(z,\xi):=z v^\prime(z) \widehat{w}(\rho)$,
and $\widehat{w}$ solves  \reff{e.what} with constants chosen so that
$$
c^*_1 I_d \ge \sigma \sigma^T, \qquad
c^*_2 I_d \ge \bar{\alpha} \bar{\alpha}^T,
\qquad
\widehat{\lambda}^i=
\tilde{\lambda}^i= 2 \overline \lambda:=
\underset{0\leq i,j\leq d}{\max}\lambda^{i,j}.
$$
Note that $W$ is explicit and is twice continuously
differentiable. We continue by showing that for 
large $K$, $V^{\eps,K}$ is a subsolution of the 
dynamic programming equation (\ref{e.dpp}), which would imply that $V^{\eps,K}\le v^\eps$ by the comparison result for the equation \reff{e.dpp} (see \cite{ks2009} Theorem 4.3.2, Proposition 4.3.4 and the subsequent discussion). A straightforward calculation, as in the proof
of Lemma 7.2 in \cite{st}, shows that the gradient
constraint
$$
\Lambda^\eps_{i,0}\cdot(V^{\eps,K}_x, V^{\eps,K}_{y_i}) \le 0,
\quad
{\mbox{ holds whenever}} 
\quad
2 \overline \lambda
+ \frac{\partial \widehat{w}}{\partial \rho_i}(\rho) \le 0.
$$
Similarly, 
$$
\Lambda^\eps_{0,i}\cdot(V^{\eps,K}_x, V^{\eps,K}_{y_i}) \le 0,
\quad
{\mbox{holds whenever}} 
\quad
-2 \overline \lambda
+ \frac{\partial \widehat{w}}{\partial \rho_i}(\rho) \le 0.
$$

Assume therefore that the elliptic
part in the equation \reff{e.what} holds.  We claim
that for a large constant $K$,
$$
\Ic(V^{\eps,K}):=\beta V^{\eps,K}
- \Lc V^{\eps,K} - \tilde{U}(V^{\eps,K}_x) \le 0.
$$
We proceed exactly as in subsection \ref{remainder} below
to arrive at
$$
\Ic(V^{\eps,K}) = \eps^2\left[-\frac12\abs{\sigma(s)\xi}v_{zz}+\frac12\Tr{\alpha\alpha^T(s,z)
W_{\xi\xi}(s,z,\xi)}-K\mathcal Au(z)+\mathcal R^\eps(z,\xi)\right].
$$
It is clear that
$\left|\mathcal R^\eps(z,\xi)\right| \le k^* \eps zv^\prime(z)$
for some constant $k^*$. We now use the elliptic part of the equation
\reff{e.what} and the choices of $c^*_i$ to conclude that
in this region,
$$
-\frac12\abs{\sigma(s)\xi}v_{zz}+\frac12\Tr{\alpha\alpha^T(s,z)
W_{\xi\xi}(s,z,\xi)} \le z v^\prime(z) \widehat{a}.
$$
Also by homotheticity, $a(z) = zv^\prime(z) a_0$ for some 
$a_0>0$.  Hence,
$$
\Ic(V^{\eps,K}) \le \eps^2 zv^\prime(z) \left[ \widehat{a} - Ka_0 + \eps k^*
\right] \le 0,
$$
provided that $K$ is sufficiently large.
Since $V^{\eps,K}$
is a smooth subsolution of the 
dynamic programming equation,
we conclude by using the standard verification argument.
\qed 

\subsubsection{Assumptions on the corrector equations}

Let $b$ be as in \reff{e.bound}, and set
\begin{equation}
\label{e.B}
B(s,z):= b\big(s,z-\ybf(s,z), \ybf(s,z)\big),
\qquad
s\in(0,\infty)^d,~z \ge 0.
\end{equation}

\begin{Assumption}[Second corrector equation: comparison]
\label{a.compare}
For any upper-semicontinuous {\rm{(}}resp. lower-semicontinuous{\rm{)}} viscosity subsolution 
{\rm{(}}resp.~supersolution{\rm{)}}
$u_1$ {\rm{(}}resp. $u_2${\rm{)}} of \reff{eq:corrector2} in $(0,\infty)^{d+1}$ satisfying the growth condition
$|u_i| \le B$ on $(0,\infty)^{d+1}$, $i=1,2$, we have $u_1 \le u_2$ in $(0,\infty)^{d+1}$.
\end{Assumption}

In the above comparison, notice that the growth of the supersolution and the subsolution is controlled by the function $B$ which is defined in \reff{e.B} by means of the local bound function $b$. In particular, $B$ controls the growth both at infinity and near the origin. 

\vspace{0.35em}

Our next assumption concerns the continuity of the solution $(w,a)$ of the first corrector equation in the parameters $(s,z)$. Recall that $w=\eta v_z\overline w$ and $a=\eta v_z\overline a$.

\begin{Assumption}[First corrector equation: regular dependence on the parameters]
\label{assump.3}
The set $\mathcal O_0$ and $a(s,z)$ are continuous in 
$(s,z)$. Moreover, $w$ is $C^2$ in $(s,z)$, and satisfies the following estimates
\begin{align}\label{estim1}
&\left(\abs{w_s}+\abs{w_{ss}}+\abs{w_z}
+\abs{w_{sz}}+\abs{w_{zz}}\right)(s,z,\xi)\leq C(s,z)\left(1+\abs{\xi}\right)\\
&\left(\abs{w_\xi}+\abs{w_{s\xi}}+\abs{w_z\xi}\right)(s,z,\xi)\leq C(s,z),
\end{align}
where $C(s,z)$ is a continuous function 
depending on the Merton value function and its derivatives.
\end{Assumption}

Notice that by the stability of viscosity solutions, the function $w$ is clearly continuous in $(s,z)$. Moreover, Assumption \ref{assump.3} is satisfied when we consider the constant coefficients and the power utility function. Indeed, in that case, there is no dependence in the $s$ variable, as emphasized in Lemma \ref{l.bound}, and the dependence in $z$ can be factored out by homogeneity.

\subsection{The first order expansion result}

The main result of this paper is the following $d-$dimensional extension of \cite{st}.

\begin{Theorem}
\label{t.main}
Under Assumptions \ref{assump2}, \ref{assump1}, 
\ref{a.compare}, and \ref{assump.3},
the sequence $\{\overline{u}^\ep\}_{\eps>0}$, defined in \reff{ubareps}, converges locally 
uniformly to the function $u$ defined in \reff{eq:corrector2}.
\end{Theorem}

\proof
In Section \ref{sect:convergence}, we will show that, the semi-limits
$u_*$ and $u^*$ are viscosity supersolution and subsolution, 
respectively, of \reff{eq:corrector2}.
Then, by the comparison Assumption \ref{a.compare}, we 
conclude that $u^* \le u \le  u_*$. Since the opposite inequality is obvious, this implies that $u^*=u_*=u$.
The local uniform convergence follows immediately
from this and the definitions.
\qed

\section{Numerical results {\normalsize (by Lo\"{i}c Richier and Bertrand Rondepierre)}}
\label{sect:numerics}

\subsection{A two-dimensional simplified model}

In this section, we report some examples of numerical results. We follow the numerical scheme suggested in Campillo \cite{campillo} which combines the finite differences approximation and the policy iteration method in order to produce an approximation of the solution of the first corrector equation. We recall that a first order approximation of the no-transaction region, for fixed $s\in(0,\infty)^d$, is deduced from the free boundary set ${\mathcal O}_0$ of Proposition \ref{prop} by:
 \be\label{hatNTeps}
 \widehat{\mbox{\bf NT}}_\eps(s)
 &:=&
 \big\{(x,y)\in\R^{d+1}: y=\theta(s,z)+\eps\eta(s,z)\mathcal{O}_0
 \big\},
 \ee
where we again denoted $z:=x+\sum_{i\le d}y^i$.
Also, as explained in Remarks \ref{rem.rep} and \ref{rem.rep2} below, the first corrector equation corresponds to a singular ergodic control problem whose corresponding optimal control processes can be viewed as a first order approximation of the optimal transfers in the original problem of optimal investment and consumption under transaction costs.

\vspace{0.35em}
For simplicity, the present numerical results are obtained in the two-dimensional case $d=2$ under the following simplifications:
 \b*
 U(c)=\frac{c^p}{p},~p\in(0,1),
 &\mu\in\R^2~~\mbox{and}~~\sigma\in \Mc_2(\R)&
 \mbox{are constants.}
 \e*
Under this simplification, it is well-known that the value function of the Merton problem $v$ is homogeneous in $z$, which in turn implies that the optimal control $\theta(z)$ is actually linear in $z$. Hence, the coefficient $\overline{\alpha}$ reduces to a constant and the solution of the first corrector equation is independent of $z$.

\subsection{Overview of the numerical scheme}
As explained above, we use the formal link established between the first corrector equation and the ergodic control problem to perform our numerical scheme. The first step in the numerical approximation for the first corrector equation is to transform the original ergodic control problem to a control problem for a Markov process in continuous time and finite state space. To do so, we restrict the domain of the variable $\rho$ to a bounded (and large enough) subset of $\mathbb R^2$, denoted by $D$, which is then discretized with a regular grid containing $N^2$ points. Then, the diffusion part of the first corrector equation is approximated using a well chosen finite difference scheme (we refer the reader to \cite{campillo} for more details). The discretized HJB equation thus takes the following form:
\b*
\min_{m\geq 0}\left(\sum_{\rho' \in \mathcal{D}}{\mathcal{L}^m_h(\rho,\rho')\overline{w}(\rho')+f^m(\rho)}\right)=\overline{a},
&\mbox{for all}&
\rho\in\mathcal{D},
\e*
where $\mathcal L^m_h$ is the discretized version of the infinitesimal generator appearing in the first order corrector equation, $m$ corresponds to the control (that is to say that $m$ is a $3\times 3$ matrix), and
\b*
f^m(\rho):=\frac{\mid\sigma^T\rho \mid^2}{2}+\Tr{\lambda^Tm},
&\lambda=(\lambda^{i,j})_{0\le i,j\le 2},&
\rho\in\mathcal{D}.
\e*
The policy iteration algorithm corresponds now to the following iterative procedure which starts from an arbitrary initial policy $m^0$, and involves the two following steps:

\vspace{0.35em}
\rm{(i)} Given a policy $m^j$, we compute for all $\rho\in D$ the solution $(\overline{w}^j,\overline{a}^j)$ of the linear system
	$$
	\sum_{\rho' \in \mathcal{D}}{\mathcal{L}^{m^j}_h(\rho,\rho')\overline{w}^j(\rho')+f^m(\rho)}=\overline{a}^j.
	$$
Notice that this is a linear system of $N^2$ equations for the $N^2$ unknowns $\big(\overline{w}^j(\rho),\rho\in\mathcal{D}\setminus\{0\}\big)$ and $\overline{a}^j$. Recall that $\overline{w}^j(0)=0$ is given.

\vspace{0.35em}
\rm{(ii)} We update the optimal control by solving the $N^2$ minimization problems
$$
m^{j+1}(\rho) \in \mbox{arg}\min_{m\geq 0}{\sum_{\rho' \in \mathcal{D}}\big\{\mathcal{L}^m_h(\rho,\rho')\overline{w}^j(\rho')+f^m(\rho)\big\}}.
$$

Finally, the algorithm is stopped whenever the difference between $\overline{a}^{j+1}$ and $\overline{a}^j$ is smaller than some initially fixed error. 

\subsection{The numerical results}

We present below different figures obtained with the above numerical scheme, representing the domain $\mathcal O_0$ and partitioning the state space in different regions depending on which controls are active or not. We remind once more the reader that $\mathcal O_0$ provides a first order approximation \reff{hatNTeps} of the true no-transaction region. For later reference, we provide our color correspondence, where we remind the reader that $0$ corresponds to the cash account and $1$ and $2$ to the two risky assets.

\begin{center}
\begin{tabular}{|c|c|c|c|c|c|c|c|c|}
\hline
{\bf  Transactions} & \tiny NT & \tiny $1$/$2$ &\tiny  $0$/$1$ & \tiny $0$/$2$ & \tiny $0$/$1$ and $0$/$2$ & \tiny $1$/$2$ and $0$/$1$ & \tiny $1$/$2$ and $0$/$2$ & \tiny $1$/$2$ and $0$/$1$ and $0$/$2$\\ \hline
{\bf Color} & &\cellcolor{blue} & \cellcolor{red} & \cellcolor{yellow}& \cellcolor{orange} & \cellcolor{violet} & \cellcolor{green} & \cellcolor{black}\\ \hline
\end{tabular}
\end{center}
and $i/j$ indicates that a positive transfer from asset $i$ to $j$, or from $j$ to $i$, occurs in the corresponding region.

\subsubsection{Cash-to-asset only}

In all existing literature addressing either numerical procedures or asymptotic expansions for the multidimensional transaction costs problem, the transactions are only allowed between a given asset and the bank account, which in our setting translates into $\lambda^{i,j}<+\infty$, if and only if $i=0$ or $j=0$, see for instance Muthuraman and Kumar \cite{muth} or Bichuch and Shreve \cite{bs2011}. In this first section, we restrict ourselves to this case and show that our numerical procedure reproduces the earlier findings. 

\vspace{0.35em}
First, we consider the following symmetric transaction costs structure with the following values for $\sigma$
$$
\lambda_0=\left(\begin{matrix} 0 & 0.001 & 0.001 \\ 0.001 & 0 & \infty \\ 0.001 & \infty & 0 \end{matrix}\right),\ 
\sigma_0=I_d,\ 
\sigma_-=\left(
  \begin{array}{ c c}
     1 & -0.25 \\
     -0.25 & 1
  \end{array} \right),\ 
\sigma_+=\left(
  \begin{array}{ c c}
     1 & 0.25 \\
     0.25 & 1
  \end{array} \right).$$

\begin{figure}[ht!]
\begin{center}
\includegraphics[scale=0.4]{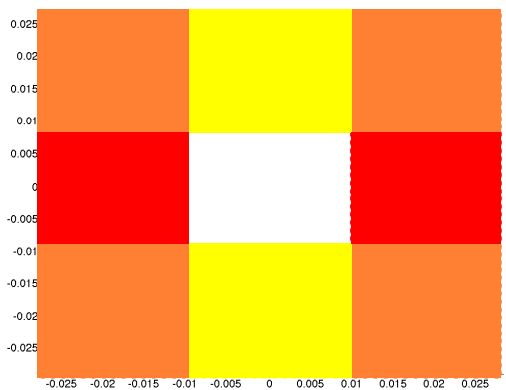}
\caption{\label{fig:cashtoassets0corr} Uncorrelated case.}
\end{center}
\end{figure}

\begin{figure}[ht!]
\begin{center}
\includegraphics[scale=0.38]{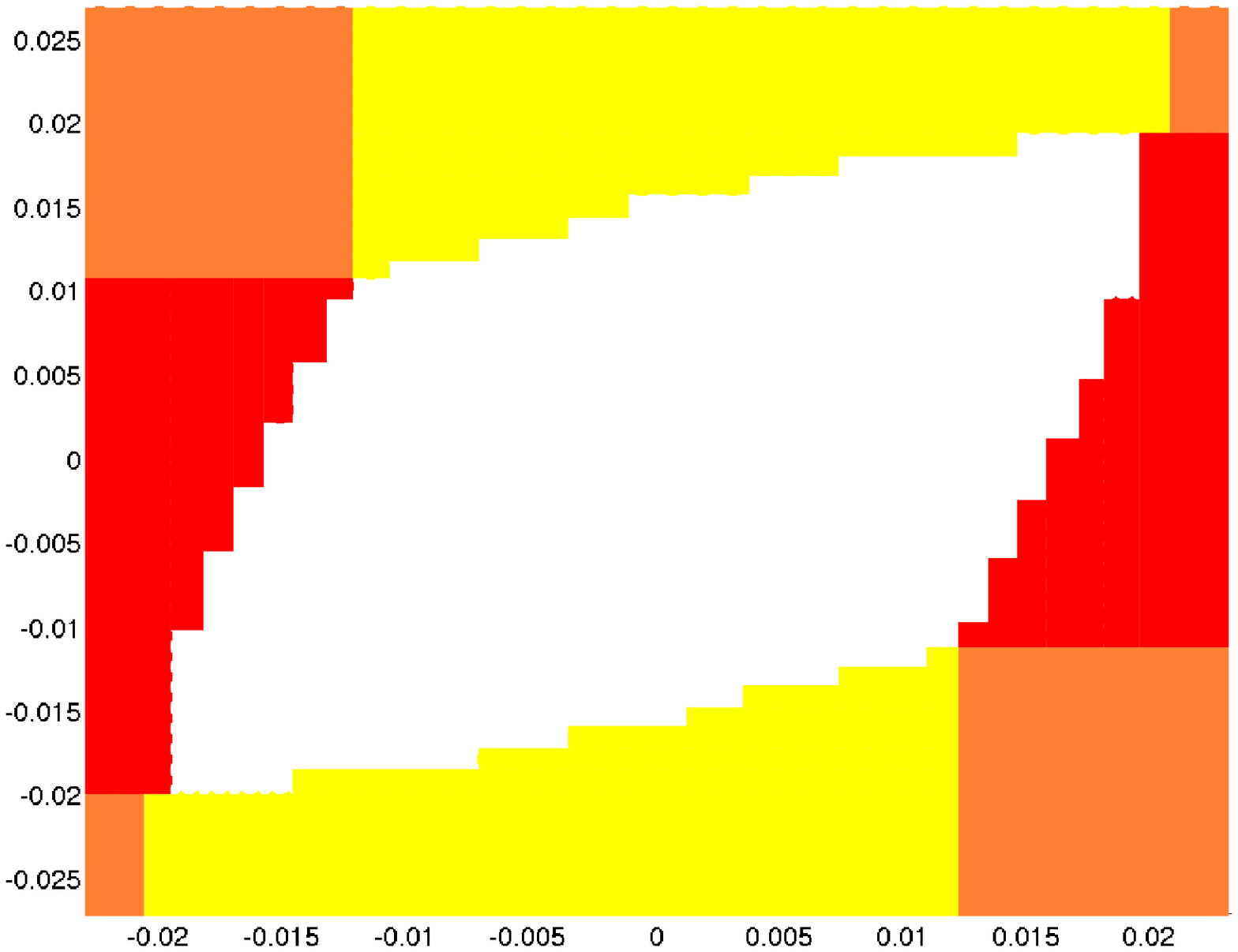}
\includegraphics[scale=0.38]{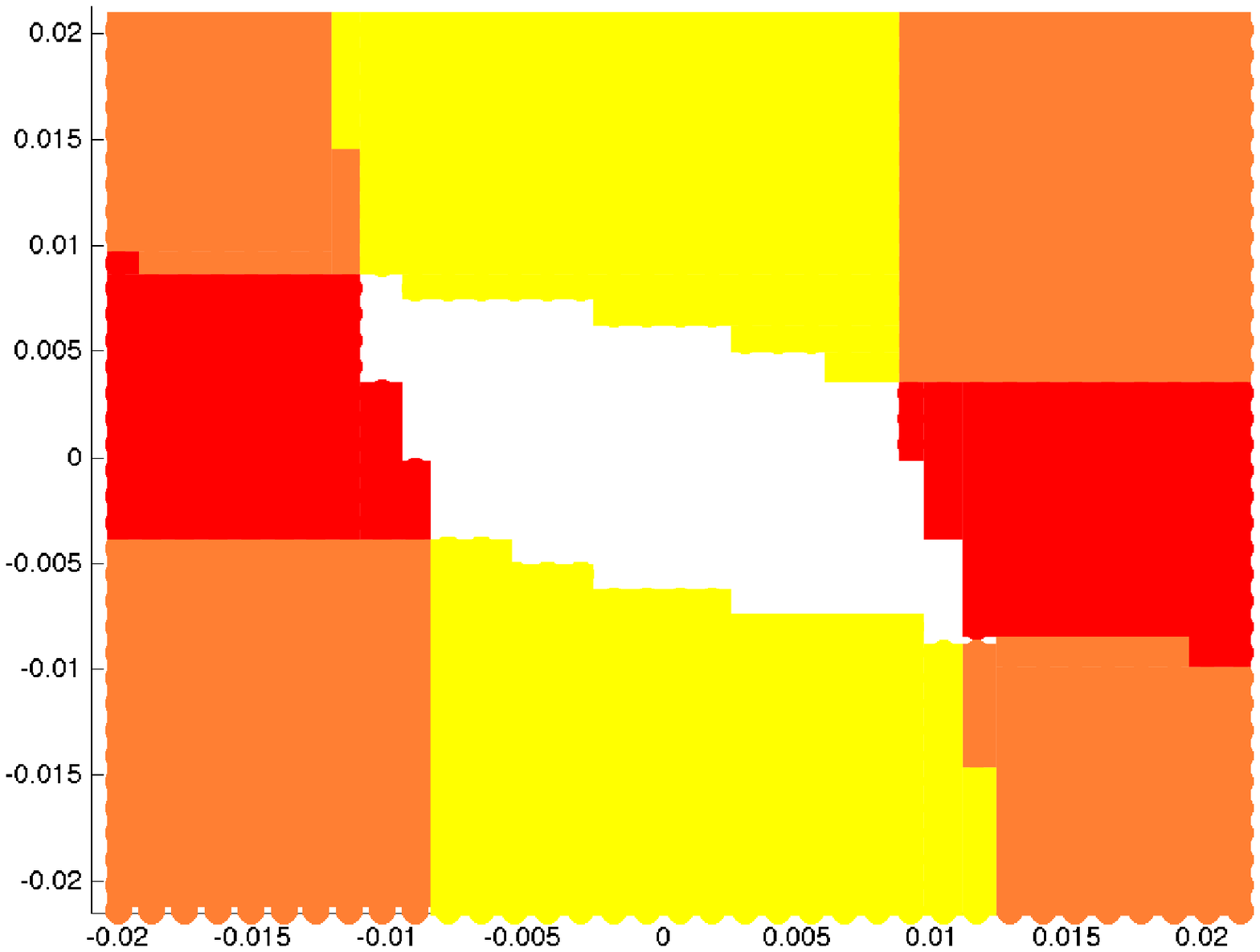}
\caption{\label{fig:cashtoassetscorr} Negative correlation (left), positive correlation (right).}
\end{center}
\end{figure}

As expected and in line with the results of \cite{muth}, the no-transaction region in the uncorrelated case is a rectangle. Then, under a possible correlation between the assets, the shape of the region is modified to a parallelogram, the direction of the deformation depending on the sign of the correlation.

\vspace{0.35em}
In the next Figure \ref{fig:highercorr}, we keep the same transaction costs structure and we consider the higher correlations induced by the volatility matrices:
 \b*
 \sigma_{--}
 =
 \left(\begin{matrix} 1 & -0.25 \\ -0.1 & 1 \end{matrix}\right) 
 &\mbox{ and }&
 \sigma_{++}
 =\left(\begin{matrix} 1 & 0.25 \\ 0.1 & 1\end{matrix}\right).
 \e*
  
  \begin{figure}[ht!]
  \begin{center}
\includegraphics[scale=0.38]{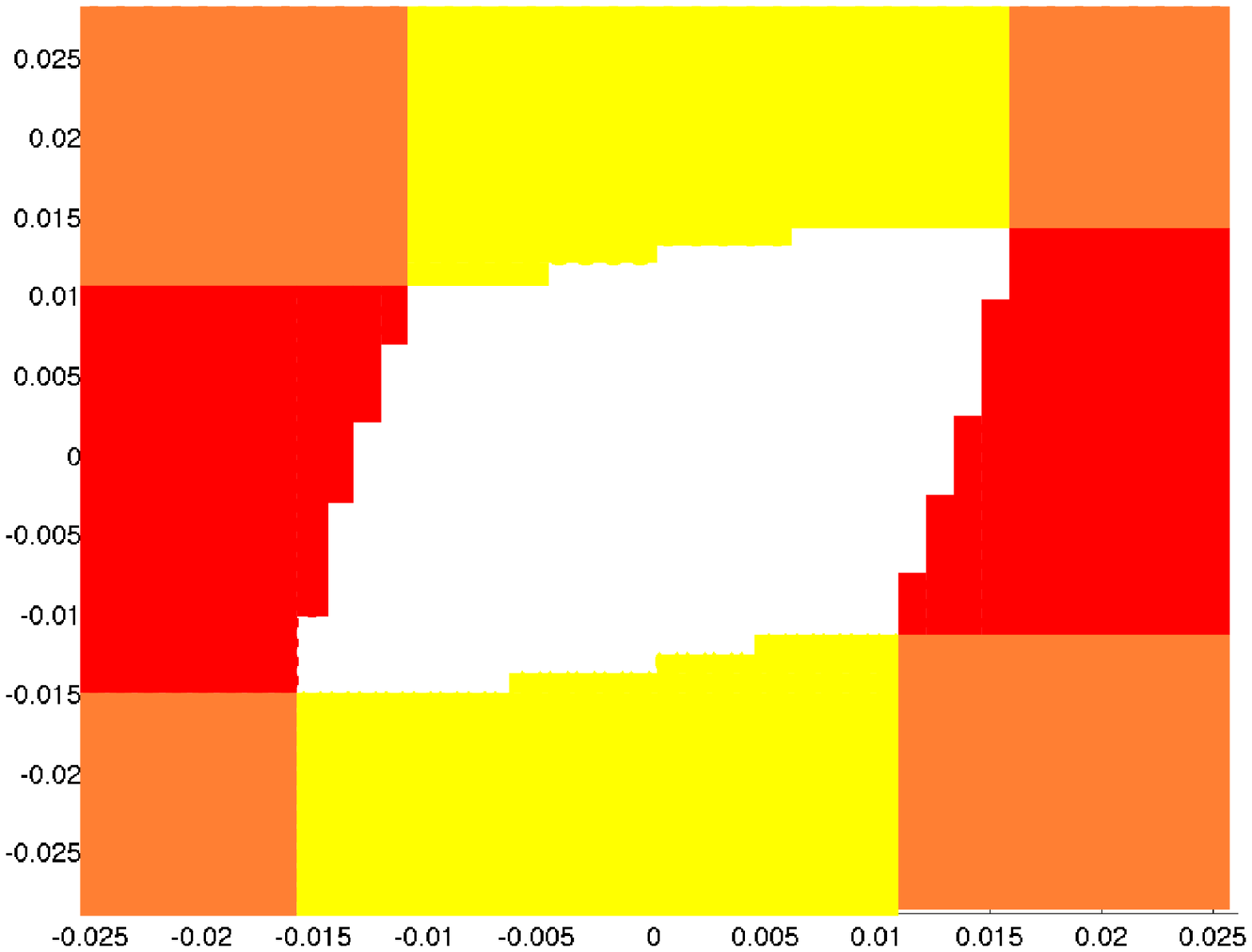}
\includegraphics[scale=0.38]{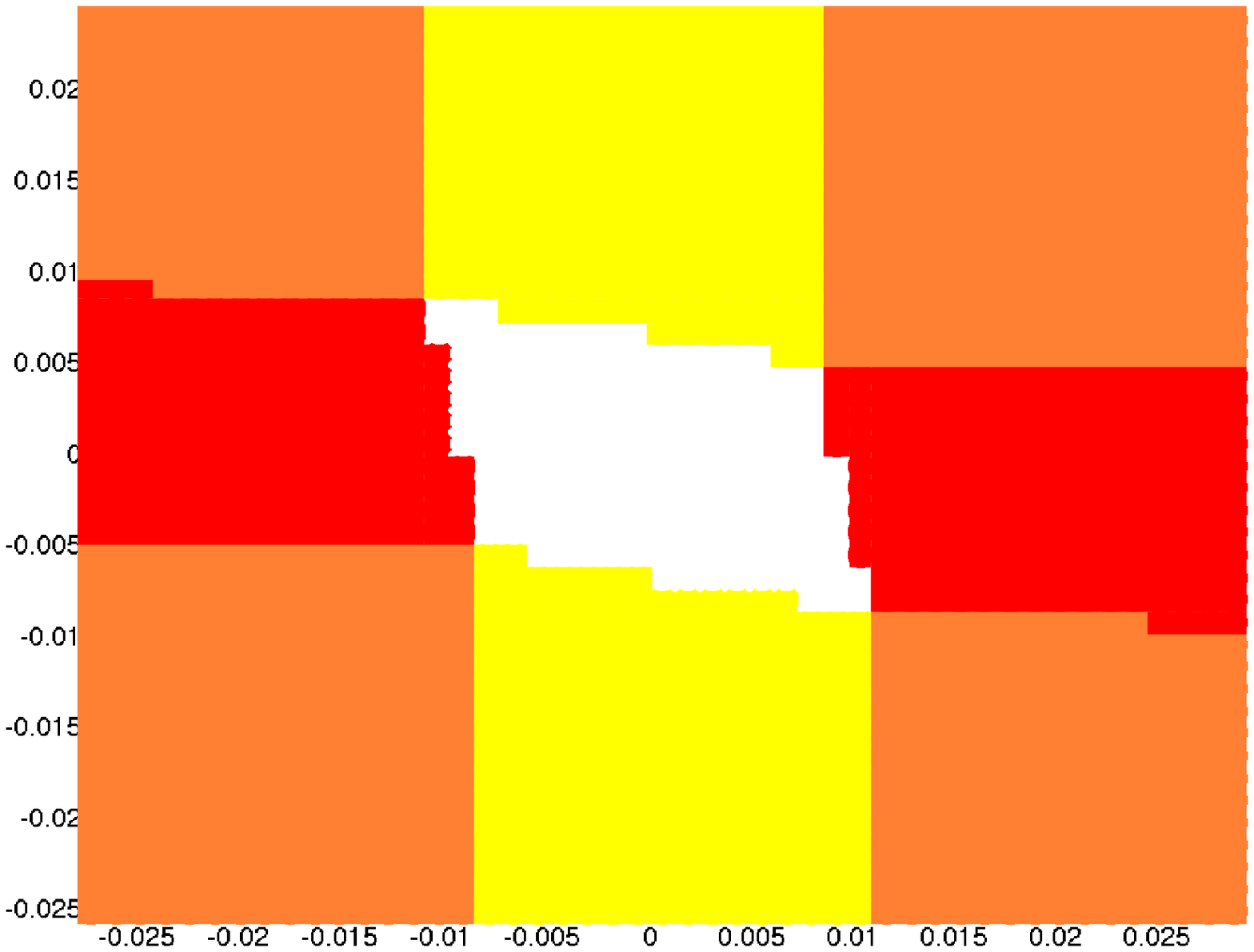}
\caption{\label{fig:highercorr} Higher correlations: negative (left) and positive (right).}
\end{center}
\end{figure}

By comparing them to the first figures, we observe that modifying the correlation induces a rotation of the no-transaction regions.

\vspace{0.35em}

\begin{figure}[ht!]
  \begin{center}
\includegraphics[scale=0.4]{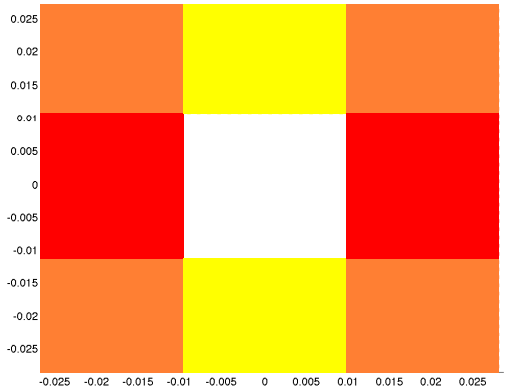}
\caption{\label{fig:lambdanonsym} Different transaction costs.}
\end{center}
\end{figure}

Our last Figure \ref{fig:lambdanonsym} for this section shows the impact of letting the transaction costs from one asset to the cash be higher than for the other asset. To isolate this effect we consider again the uncorrelated case $\sigma_0$, and choose $
\lambda'=\left(\begin{matrix} 0 & 0.001 & 0.002 \\ 0.001 & 0 & \infty \\ 0.002 & \infty & 0 \end{matrix}\right).$

As expected, we still observe a rectangle but with modified dimensions. More precisely, transactions between the first asset and the cash account occur more often since they are cheaper.

\subsubsection{Possible transfers between all assets}

In this section we allow for transactions between all assets, a feature which was not considered in any of the existing numerical approximations in the literature on the present problem. We start by fixing a symmetric transaction costs structure $\lambda_0$, and we illustrate the impact of correlation by considering the volatility matrices $\sigma_0$, $\sigma_-$ and $\sigma_+$.

\begin{figure}[ht!]
\begin{center}
\includegraphics[scale=0.4]{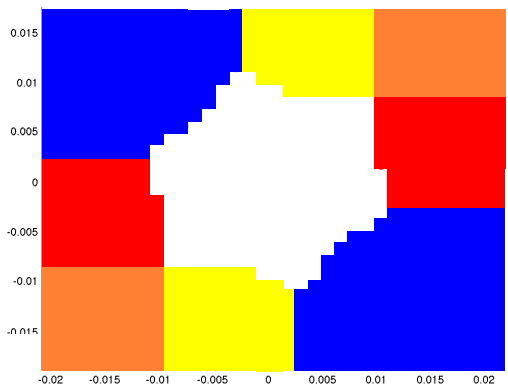}
\caption{\label{fig:allassets0corr} Uncorrelated case.}
\end{center}
\end{figure}

\begin{figure}[ht!]
\begin{center}
\includegraphics[scale=0.38]{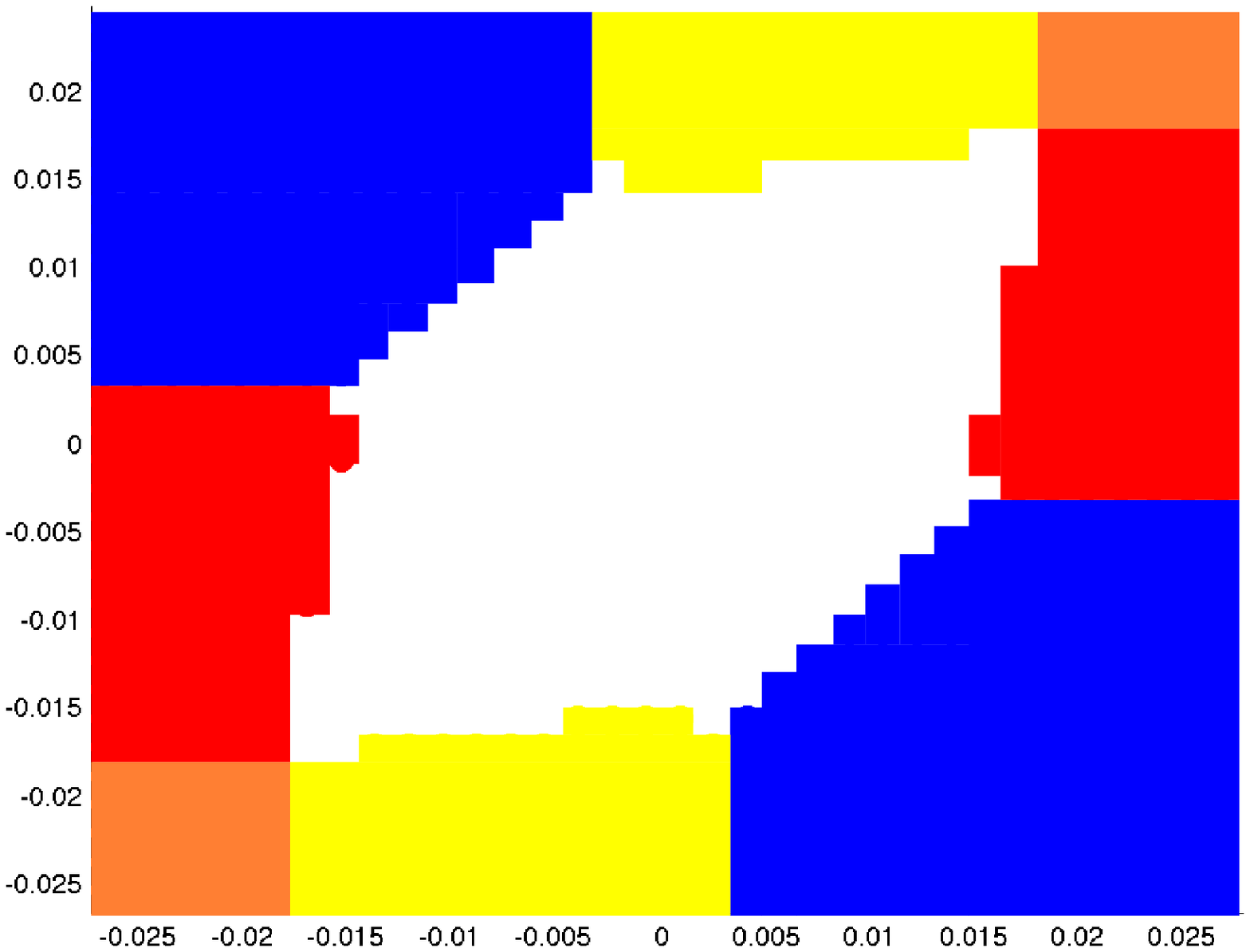}
\includegraphics[scale=0.38]{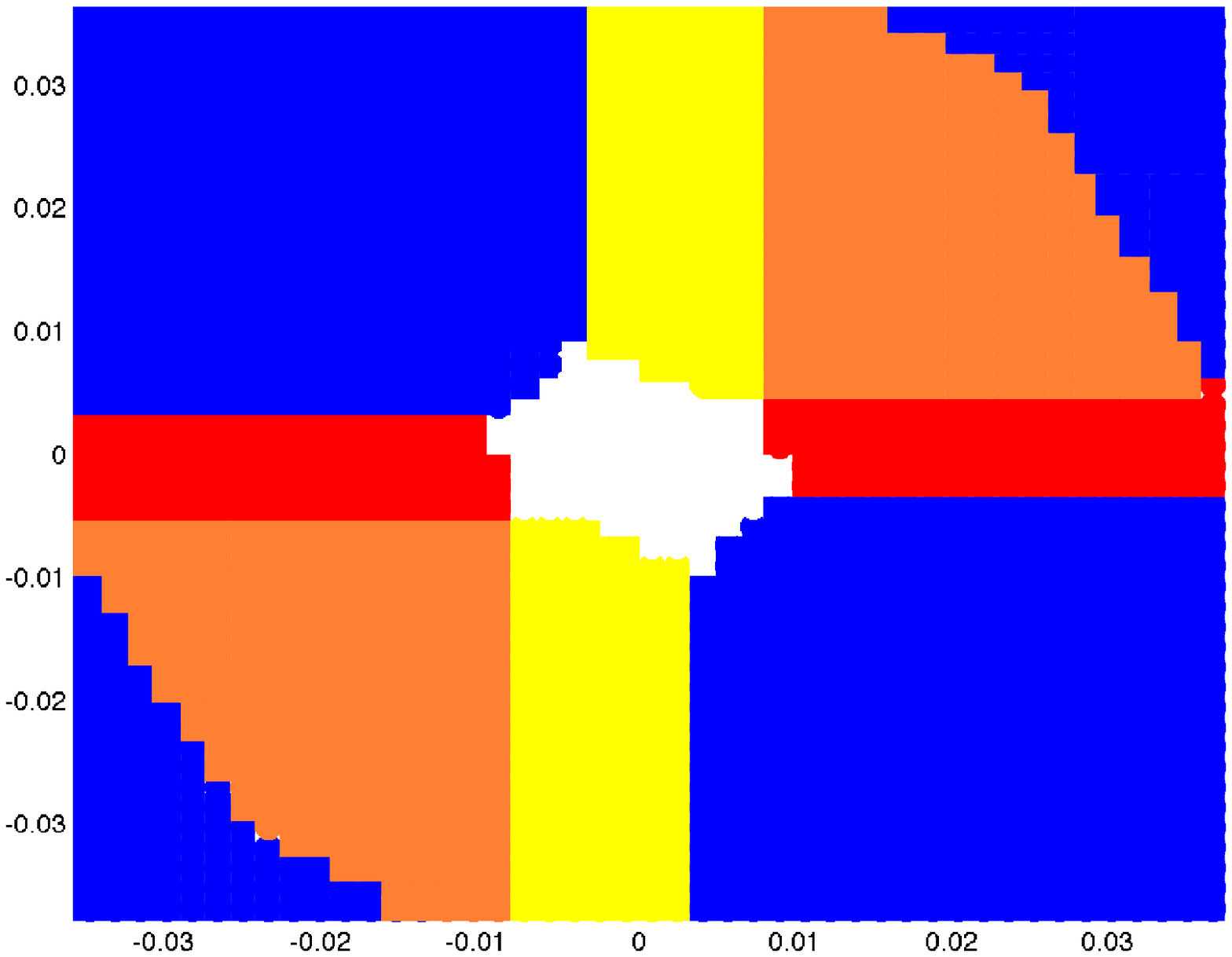}
\caption{\label{fig:allassetscorr} Negative correlation (left), positive correlation (right).}
\end{center}
\end{figure}

\vspace{0.35em}
Our first observation is that transactions between the two assets do occur, and more importantly that as a consequence, the no transaction region seems to no longer be convex, an observation which, as far as we know, was not made before in the literature. Moreover, as in the previous section, the introduction of correlation induces a deformation of the no-transaction region. In order to insist on this loss of convexity, we also report the following figure obtained with the same parameters as the left side of Figure \reff{fig:allassets0corr}, but with more precise computations and without colors

  \begin{figure}[ht!]
  \begin{center}
\includegraphics[scale=0.4]{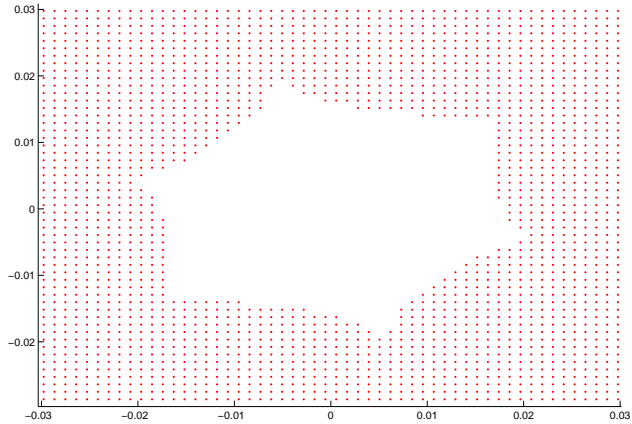}
\caption{All transactions allowed.}
\end{center}
\end{figure}

\vspace{0.35em}
Our last figure shows the impact of an asymmetric transaction costs structure. The volatility matrix is set to $\sigma_0$ and the transaction costs matrix to $\lambda'$. As expected, there are almost no transactions between the cash account and asset $1$, since they are twice as expensive as the other ones. Surprisingly, we also observe the occurrence of small zones (in black and violet), where transactions are simultaneously performed between the assets and between the assets and the cash account.
  \begin{figure}[ht!]
  \begin{center}
\includegraphics[scale=0.4]{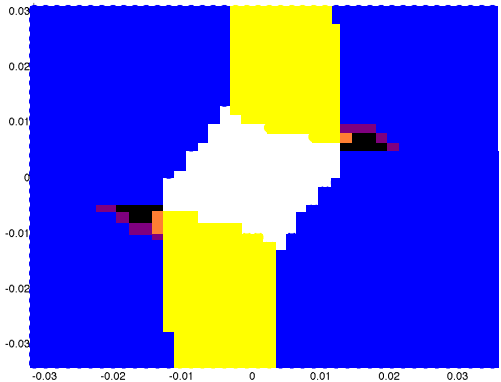}
\caption{Asymmetric transaction costs.}
\end{center}
\end{figure}

\section{Convergence}
\label{sect:convergence}

In the rest of the paper, we denote for any function $f(s,x,y)$:
 $$
 \hat f(s,z,\xi)
 :=
 f\big(s,z-\ybf(s,z)\cdot\1_d,\eps\xi+\ybf(s,z)\big).
 $$
This section is dedicated to the proof of our main result, Theorem \ref{t.main}. Let:
$$u^\eps(s,x,y):=\bar{u}^\eps(s,x,y)-\eps^2w(s,z,\xi), \ s\in\mathbb R_+^d,\ (x,y)\in K_\eps.$$

\subsection{First estimates and properties}
We start by obtaining several estimates of $u^\eps$. 
Set
$$
\overline \lambda:=
\underset{0\leq i,j\leq d}{\max}\lambda^{i,j},\ \underline \lambda
:=\underset{0\leq i,j\leq d}{\min}\lambda^{i,j}.
$$
We also recall that $L$ is the upper bound of the set $C$.

\begin{Lemma}\label{lemma.local}
For $(\eps,s,x,y)\in(0,+\infty)\times(0,+\infty)^d\times K_\eps$, and $z:=x+y$, we have
	$$u^\eps(s,x,y)\geq -\eps L v_z(s,z)\abs{y-\ybf(s,z)}.$$
Consequently, under Assumption \ref{assump1}, we have for all $(\eps,s,x,y)
\in(0,+\infty)\times(0,+\infty)^d\times K_\eps$
	$$0\leq u_*(s,x,y)\leq u^*(s,x,y)<+\infty.$$
\end{Lemma}

\proof Since $v^\eps(s,x,y)\leq v(s,z)$, it follows from the definition of $u^\eps$ that 
$$
u^\eps(s,x,y)\geq -\eps^2w(s,z,\xi).
$$
Next, recall that $D\overline{w}$ takes values in the bounded set $C$. Since $w(.,0)=0$, this implies that $-w(s,z,\xi)\geq -L\abs{\xi}v_z(s,z)$, and completes the proof.
\qed

\vspace{0.4em}

The next Lemma proves that the relaxed semi-limits are only functions of $(s,z)$.

\begin{Lemma}\label{indep}
Let Assumptions \ref{assump2} and \ref{assump1} hold. Then,
$u^*$ and $u_*$ are only functions of $(s,z)$. Furthermore, we have
 \b*
 u_*(s,z)
 &=&
 \underset{(\eps,s',z')\rightarrow (0,s,z)}{\underline{\lim}}
 u^\eps\big(s',z'-\ybf(s',z')\cdot \1_d,\ybf(s',z')\big)
 \\
 u^*(s,z)
 &=&
 \underset{(\eps,s',z')\rightarrow (0,s,z)}{\overline{\lim}}
 u^\eps\big(s',z'-\ybf(s',z')\cdot \1_d,\ybf(s',z')\big).
 \e*
\end{Lemma}

\proof
We proceed in several steps. We assume throughout the proof that the parameter $\eps$ is less than one.

\vspace{0.35em}
{\bf Step $1$:} In view of 
the gradient constraints of the dynamic programming 
equation satisfied by $v^\eps$, we know that for all
$1\leq i\leq d$, we have in the viscosity sense
\begin{equation}
\Lambda^\eps_{i,0}.(v^\eps_x,v^\eps_y)\geq 0,\ \Lambda^\eps_{0,i}.(v^\eps_x,v^\eps_y)\geq 0.
\label{eq:visco}
\end{equation}

Now define 
$$\widehat v^\eps(s,z,\xi):=v^\eps(s,z-\eps\xi\cdot \1_d-\ybf(s,z)\cdot \1_d,\eps\xi+\ybf(s,z)).$$

We directly calculate that \reff{eq:visco} implies for all 
$1\leq i\leq d$ in the viscosity sense
\begin{equation}
\eps^4\lambda^{i,0}\widehat{v}^\eps_z(s,z,\xi)-\eps^3\lambda^{i,0}\ybf_z(s,z).\widehat{v}^\eps_\xi(s,z,\xi)
+(1+\eps^3\lambda^{i,0})\widehat{v}^\eps_{\xi^i}(s,z,\xi)\geq 0,
\label{eq:1}
\end{equation}
and
\begin{equation}
\eps^4\lambda^{0,i}\widehat{v}^\eps_z(s,z,\xi)
-\eps^3\lambda^{0,i}\ybf_z(s,z).\widehat{v}^\eps_\xi(s,z,\xi)-\widehat{v}^\eps_{\xi^i}\geq 0.
\label{eq:2}
\end{equation}
Using the fact that all the $\ybf^i_z(s,z)$ are strictly positive 
(see Assumption \ref{assump2}), we can multiply \reff{eq:1} 
by $\ybf^i_z(s,z)$ and sum to obtain, once more in the viscosity sense
\begin{equation}
\left(1-\eps^3 \sum_{i=1}^d\frac{\ybf^i_z(s,z)\lambda^{i,0}}
{1+\eps^3\lambda^{i,0}}\right)\ybf_z(s,z).\widehat{v}^\eps_\xi(s,z,\xi)
\geq -\eps^4\sum_{i=1}^d
\frac{\lambda^{i,0}}{1+\eps^3\lambda^{i,0}}\widehat{v}^\eps_z(s,z,\xi).
\label{eq:3}
\end{equation}
Now, we have by Assumption \ref{assump2}
$$
1-\eps^3 \sum_{i=1}^d\frac{\ybf^i_z(s,z)\lambda^{i,0}}{1+\eps^3\lambda^{i,0}}
=1-\sum_{i=1}^d\ybf_z^i(s,z)+\sum_{i=1}^d\frac{\ybf^i_z(s,z)}{1+\eps^3\lambda^{i,0}}\geq 0.
$$
Using this inequality in \reff{eq:3} yields, in the viscosity sense
\begin{equation}
\ybf_z(s,z).\widehat{v}^\eps_\xi(s,z,\xi)\geq 
- \frac{\sum_{i=1}^d\frac{\lambda^{i,0}}{1+\eps^3\lambda^{i,0}}}{1-\eps^3
 \sum_{i=1}^d\frac{\ybf^i_z(s,z)\lambda^{i,0}}{1+\eps^3\lambda^{i,0}}}\eps^4\widehat{v}^\eps_z(s,z,\xi).
\end{equation}
Plugging this estimate in \reff{eq:1} and \reff{eq:2}, we obtain in the viscosity sense
$$
\widehat v^\eps_{\xi^i}(s,z,\xi)\leq \lambda^{0,i}\eps^4
\left(1+\frac{\sum_{i=1}^d\frac{\lambda^{i,0}}{1+\eps^3\lambda^{i,0}}}{1-\eps^3 
\sum_{i=1}^d\frac{\ybf^i_z(s,z)\lambda^{i,0}}{1+\eps^3\lambda^{i,0}}}\right)\widehat v^\eps_z(s,z,\xi)
\leq \eps^4\overline\lambda\left(1+\frac{\overline{\lambda}d}{c_0}\right)\widehat v^\eps_z(s,z,\xi),
$$
and
$$\widehat v^\eps_{\xi^i}(s,z,\xi)\geq -2\eps^4\underline{\lambda}\widehat v^\eps_z(s,z,\xi).$$

By the concavity of $v^\eps$, its gradient exists almost everywhere.
 Moreover, since we assumed that $\ybf$ is smooth, 
this implies that $\widehat v^\eps_z$ also exists almost everywhere. 
Hence, we conclude from the above estimates that
 \be\label{eq:5}
 \abs{\widehat v^\eps_\xi}\leq A\eps^4\widehat v^\eps_z,
 &\mbox{where}&
 A
 :=
 \sqrt{d}\max\Big\{\overline\lambda
                   \Big(1+\frac{\overline{\lambda}^2d}{c_0}\Big),
                   2\underline\lambda
             \Big\}.
 \ee
{\bf Step $2$:} We now prove an estimate for $\widehat v^\eps_z$. 
By definition and from Assumption \ref{assump2}, 
\begin{align}\label{eq:8}
\nonumber\widehat v^\eps_z(s,z,\xi)&=\partial_zv^\eps\left(s,z-\eps\xi\cdot \1_d-\ybf(s,z)\cdot \1_d,\eps\xi+\ybf(s,z)\right)\\
\nonumber&=\left(1-\ybf_z(s,z)\cdot \1_d\right)v^\eps_x(s,x,y)+\ybf_z(s,z).v^\eps_y(s,x,y)\\
&\leq v^\eps_x(s,x,y)+v_y^\eps(s,x,y)\cdot \1_d.
\end{align}

Therefore, we can focus on obtaining estimates on $v^\eps_x$ and $v^\eps_y$. First, by concavity of $v^\eps$ in $x$ and of $v$ in $z$, we have
\begin{align*}
v^\eps_x(s,x,y)&\leq \frac{v^\eps(s,x,y)-v^\eps(s,x-\eps,y)}{\eps}\\
&\leq \frac{v(s,z)-v(s,z-\eps)}{\eps}+ \frac{v(s,z-\eps)-v^\eps(s,x-\eps,y)}{\eps}\\
&\leq v_z(s,z-\eps)+\frac{v(s,z-\eps)-v^\eps(s,x-\eps,y)}{\eps}.
\end{align*}

Now, using the definition of $u^\eps$, we obtain
$$v^\eps_x(s,x,y)\leq v_z(s,z-\eps)+\eps\left(u^\eps(s,x-\eps,y)+\eps^2w(s,z-\eps,\xi_\eps)\right),$$
where
$$\xi_\eps:=\frac{y-\ybf(s,z-\eps)}{\eps}=\xi+ \frac{\ybf(s,z)-\ybf(s,z-\eps)}{\eps}.$$

From the estimates on $w$ in Theorem \ref{prop.uni} and Assumption \ref{assump.3}, we have
$$\abs{w(s,z-\eps,\xi_\eps)}\leq Lv_z(s,z)(1+\abs{\xi_\eps})\leq Lv_z(s,z)\left(1+\abs{\xi}+\frac{\abs{\ybf(s,z)-\ybf(s,z-\eps)}}{\eps}\right),$$
and therefore
 $$v^\eps_x(s,x,y)\leq v_z(s,z-\eps)+\eps u^\eps(s,x-\eps,y)+\eps^3Lv_z(s,z)\left(1+\abs{\xi}+\frac{\abs{\ybf(s,z)-\ybf(s,z-\eps)}}{\eps}\right).
 $$
As for $v^\eps_y$, we use again the concavity of $v^\eps$ and $v$ in $y$ and $z$, respectively,
\begin{align*}
v^\eps_{y^i}(s,x,y)&\leq \frac{v^\eps(s,x,y)-v^\eps(s,x,y-\eps e_i)}{\eps}\\
&\leq \frac{v(s,z)-v(s,z-\eps)}{\eps}+ \frac{v(s,z-\eps)-v^\eps(s,x,y-\eps e_i)}{\eps}\\
&\leq v_z(s,z-\eps)+\frac{v(s,z-\eps)-v^\eps(s,x,y-\eps e_i)}{\eps},
~~0\le i\le d.
\end{align*}

Similarly as above, this yields
 \b*
 v^\eps_{y^i}(s,x,y)
 &\leq& 
 v_z(s,z-\eps)+\eps u^\eps(s,x,y-\eps e_i)
 \\
 &&
 +\eps^3Lv_z(s,z)\left(1+\abs{\xi}+\frac{\abs{-\eps e_i+\ybf(s,z)-\ybf(s,z-\eps)}}{\eps}\right).
 \e*
Plugging the above estimates in \reff{eq:8}, we get
 \begin{align}\label{eq:9}
 \widehat v^\eps_z(s,z,\eps)\leq& (1+d)v_z(s,z-\eps) +\eps\left(u^\eps(s,x-\eps,y)+\sum_{i=1}^du^\eps(s,x,y-\eps e_i)\right)
\\
&+\eps^3Lv_z(s,z)\left(d+(1+d)\left(1+\abs{\xi}+\frac{\abs{\ybf(s,z)-\ybf(s,z-\eps)}}{\eps}\right)\right)=:\gamma^\eps(s,z,\xi).
\nonumber
\end{align}

{\bf Step $3$} Recall that $\widehat v^\eps_z$ exists almost everywhere, and 
thus by definition of $u^\eps$, this also holds for $\widehat u^\eps_z$ and 
$\widehat u^\eps_\xi$. Now, using \reff{eq:5}, \reff{eq:9} together with 
Assumption \ref{assump2} and our estimates on $w$, we obtain for some constant $C>0$
$$
\abs{\widehat u^\eps_\xi(s,z,\xi)}\leq \eps^2C\left(v_z(s,z)
+\widehat v^\eps_z(s,z,\xi)\right)\leq \eps^2C\left(v_z(s,z)+\gamma^\eps_z(s,z,\xi)\right).
$$

With this estimate, we can conclude the proof exactly as in the proof of Lemma $6.2$ in \cite{st}.
\qed

\subsection{Remainder estimate}\label{remainder}

In this section, we isolate an estimate which will be needed at various occasions in the subsequent proofs. The following calculation extends the estimate of Section 4.2 in \cite{st}. 
For any function 
$$
\Psi^\eps(s,x,y):=v(s,z)-\eps^2\phi(s,z)-\eps^4\varpi(s,z,\xi),
$$ 
with smooth $\phi$ and $\upsilon$ such that $\upsilon$ also verifies the estimates \reff{estim1}, we have
\begin{align*}
\mathcal I(\Psi^\eps)(s,x,y):&=\left(\beta\Psi^\eps-\mathcal L\Psi^\eps-\widetilde U(\Psi^\eps_x)\right)(s,x,y)\\
&=\eps^2\left[-\frac12\abs{\sigma(s)\xi}v_{zz}+\frac12\Tr{\alpha\alpha^T(s,z)
\varpi_{\xi\xi}(s,z,\xi)}-\mathcal A\phi(s,z)+\mathcal R^\eps(s,z,\xi)\right].
\end{align*}
Similarly as in \cite{st}, direct but tedious calculations provides the following estimate:
\begin{align*}
\abs{\mathcal R^\eps(s,x,y)}\leq 
&\eps\Big(\abs{\mu-r\cdot \1_d}\abs{\xi}\abs{\phi_z}
  +\frac{\abs{\sigma}^2}{2}
    \left(2\abs{\ybf}\abs{\xi}+\eps\xi^2\right)
    \abs{\phi_{zz}}
  +\abs{\sigma}^2\abs{\xi}\abs{D_{sz}\phi}
  \Big)(s,z)\\
&+\eps C(s,z)\left(1+\eps\abs{\xi}+\eps^2\abs{\xi}^2+\eps^3\abs{\xi}^3\right)\\
&+\eps^{-2}\abs{\widetilde{U}(\Psi^\eps_x)-\widetilde U(v_z)-(\Psi^\eps_x-v_z)\widetilde{U}^{'}(v_z)},
\end{align*}
for some continuous function $C(s,z)$. Now using the fact $\widetilde U$ is $C^1$ and convex and the estimates assumed for $\upsilon$, we obtain
 \begin{align*}
 \abs{\mathcal R^\eps(s,x,y)}
 \leq 
 &\eps\Big(\abs{\mu-r\cdot \1_d}\abs{\xi}\abs{\phi_z}
            +\frac{\abs{\sigma}^2}{2}
                   \left(2\abs{\ybf}\abs{\xi}+\eps\xi^2\right)
                   \abs{\phi_{zz}}
            +\abs{\sigma}^2\abs{\xi}\abs{D_{sz}\phi}
      \Big)(s,z)
 \\
 &+\eps C(s,z)\left(1+\eps\abs{\xi}+\eps^2\abs{\xi}^2+\eps^3\abs{\xi}^3\right)\\
&+\eps^{2}\left(\abs{\phi_z}+\eps C(s,z)(1+\eps\abs{\xi})\right)^2\widetilde{U}^{''}\left(v_z+\eps^2\abs{\phi_z}
+\eps^3C(s,z)(1+\eps\abs{\xi})\right).
\end{align*}

\subsection{Viscosity subsolution property}

In this Section, we prove
\begin{Proposition}
Under Assumptions \ref{assump2}, \ref{assump1} and \ref{assump.3}, 
the function $u^*$ is a viscosity subsolution of the second corrector equation.
\end{Proposition}

\proof
Let $(s_0,z_0,\phi)\in(0,+\infty)^d\times(0,+\infty)
\times C^2\left((0,+\infty)^d\times(0,+\infty)\right)$ be such that
\begin{equation}
(u^*-\phi)(s_0,z_0)>(u^*-\phi)(s,z),\text{ for all $(s,z)
\in(0,+\infty)^d\times(0,+\infty)\backslash\left\{(s_0,z_0)\right\}$}.
\label{eq:10}
\end{equation}

By definition of viscosity subsolutions, we need to show that 
$$\mathcal A\phi(s_0,z_0)-a(s_0,z_0)\leq 0.$$

We will proceed in several steps.

\vspace{0.35em}
{\bf Step $1$:} First of all, we know from Lemma \ref{indep} that 
there exists a sequence $(s^\eps,z^\eps)$ which realizes the 
$\limsup$ for $\widehat u^\eps$, that is to say
$$
(s^\eps,z^\eps)\underset{\eps\downarrow 0}{\longrightarrow} 
(s_0,z_0)\text{ and }\widehat u^\eps(s^\eps,z^\eps,0)\underset{\eps\downarrow 0}{\longrightarrow}u^*(s_0,z_0).
$$

It follows then easily that $l^\eps_*:=\widehat u^\eps(s^\eps,z^\eps,0)-\phi(s^\eps,z^\eps)\underset{\eps\downarrow 0}{\longrightarrow}0,$ and
$$(x^\eps,y^\eps)=\left(z^\eps-\ybf(s^\eps,z^\eps)\cdot \1_d,\ybf(s^\eps,z^\eps)\right)\underset{\eps\downarrow 0}
{\longrightarrow}(x_0,y_0):=\left(z_0-\ybf(s_0,z_0)\cdot \1_d,\ybf(s_0,z_0)\right).$$

Now recall from Assumption \ref{assump1} that $u^\eps$ is 
locally bounded from above. This implies the existence of $r_0:=r_0(s_0,x_0,y_0)>0$ 
and $\eps_0:=\eps_0(s_0,x_0,y_0)>0$ verifying
\begin{equation}
b_*:=\sup\left\{u^\eps(s,x,y),\ (s,x,y)\in B_0,\ \eps\in(0,\eps_0]\right\}<+\infty,
\label{eq:11}
\end{equation}
where $B_0:=B_{r_0}(s_0,x_0,y_0)$ is the open ball with radius $r_0$ and center $(s_0,x_0,y_0)$. Moreover, notice that we can always decrease $r_0$ so that $r_0\leq z_0/2$, 
which then implies that $B_0$ does not cross the line $z=0$. Now for any $(\eps,\delta)\in(0,1]^2$, we define $\Psi^{\eps,\delta}$ and the corresponding $\hat\Psi^{\eps,\delta}$ by
 $$
 \widehat{\Psi}^{\eps,\delta}(s,z,\xi)
 :=
 v(s,z)-\eps^2\left(l_*^\eps+\phi(s,z)
 +\widehat{\Phi}^\eps(s,z,\xi)\right)-\eps^4(1+\delta)w(s,z,\xi),
 $$
where the function $\hat\Phi^\eps$ and the corresponding $\Phi^\eps$ are given by:  
 $$
 \hat\Phi^\eps(s,x,y)
 :=
 c_0\Big((s-s^\eps)^4+(z-z^\eps)^4+\eps^4w^4(s,z,\xi)
    \Big),
 $$
and $c_0>0$ is a constant chosen large enough in order to have for $\eps$ small enough
\begin{equation}
\Phi^\eps\leq 1+b_*-\phi,\text{ on }B_0\backslash B_1,\text{ where }B_1:=B_{\frac{r_0}{2}}(s_0,x_0,y_0).
\label{eq:12}
\end{equation}
We emphasize that the constant $c_0$ may depend on $(\phi,s_0,x_0,y_0,\delta)$ 
but not on $\eps$, and that {\sl a priori} the function $\widehat \Psi^{\eps,\delta}$ 
is not $C^2$ in $\xi$, because the function $w$ is only in $C^{1,1}$. 
This is a major difference with the one-dimension  case treated in \cite{st}.

\vspace{0.35em}
{\bf Step $2$:} We now prove that for all sufficiently small $\eps$ and $\delta$, 
the difference $(v^\eps-\Psi^{\eps,\delta})$ has a local minimizer in $B_0$. First, notice that this is equivalent to showing that the following quantity has a local minimizer in $B_0$
\begin{align*}
I^{\eps,\delta}(s,x,y):&=\frac{v^\eps(s,x,y)-\Psi^{\eps,\delta}(s,x,y)}{\eps^2}\\
&=-u^\eps(s,x,y)+l_*^\eps+\phi(s,z)+\phi^\eps(s,x,y)+\eps^2\delta w(s,z,\xi).
\end{align*}
By \reff{eq:12} and the fact that $w\geq 0$, we have for any $(s,x,y)\in\partial B_0$
$$I^{\eps,\delta}(s,x,y)\geq -u^\eps(s,x,y)+l_*^\eps+1+b_*+\eps^2\delta w(s,z,\xi)\geq 1+l_*^\eps>0,$$
for $\eps$ small enough. Moreover, since $I^{\eps,\delta}(s^\eps,x^\eps,y^\eps)=0$, this implies that $I^{\eps,\delta}$ has a local minimizer $(\tilde{s}^\eps,\tilde{x}^\eps,\tilde{y}^\eps)$ in $B_0$, and we introduce the corresponding
$$
\tilde{z}^\eps:=\tilde{x}^\eps+\tilde{y}^\eps\cdot \1_d, \text{ and }\tilde{\xi}^\eps
:=\frac{\tilde{y}^\eps-\ybf(\tilde{s}^\eps,\tilde{z}^\eps)}{\eps}.
$$
We then have
$$
\underset{(s,z,\xi)}{\min}(\widehat v^\eps-\widehat \Psi^{\eps,\delta})(s,z,\xi)
=(\widehat v^\eps-\Psi^{\eps,\delta})(\tilde{s}^\eps,\tilde{z}^\eps,\tilde{\xi}^\eps)\leq 0,\ 
\abs{\tilde{s}^\eps-s_0}+\abs{\tilde{z}^\eps-z_0}\leq r_0,\ \abs{\tilde{\xi}^\eps}\leq \frac{r_1}{\eps},
$$
for some constant $r_1$. We now use the viscosity supersolution property of $v^\eps$. Since $\Psi^{\eps,\delta}$ is $C^1$, we obtain from the first order operator in the dynamic programming equation that:
 \be\label{eq:14}
 \Lambda^\eps_{i,j}
 \cdot \big(\Psi_x^{\eps,\delta},\Psi_y^{\eps,\delta}\big)
 (\tilde{s}^\eps,\tilde{x}^\eps,\tilde{y}^\eps)
 \geq 0
 &\mbox{for all}&
 0\leq i,j\leq d.
 \ee
{\bf Step $3$:} In this step, we show that for $\eps$ small enough, we have
\begin{equation}
\tilde{\rho}^\eps
\;:=\;
\frac{\tilde{\xi}^\eps}{\eta(\tilde{s}^\eps,\tilde{z}^\eps)}
\;\in\;
\mathcal O_0(\tilde{s}^\eps,\tilde{z}^\eps).
\label{eq:13}
\end{equation}
where $\mathcal O_0(s,z)$ is the open set of Proposition \ref{prop}. 
We argue by contradiction assuming that there exists some sequence $\eps_n\longrightarrow 0$ such that $\tilde{\rho}^{\eps_n} \not \in \mathcal O_0(\tilde{s}^{\eps_n},\tilde{z}^{\eps_n})$.
This implies that
 \be\label{eq:15}
 -\lambda^{i_0^n,j_0^n}
 +\big(\partial_{i_0^n}\overline{w}-\partial_{j_0^n}\overline{w}\big)
  (\tilde{s}^{\eps_n},\tilde{z}^{\eps_n},\tilde{\rho}^{\eps_n})=0
 &\mbox{for some}&
 (i_0^n,j_0^n).
 \ee
By the gradient constraint \reff{eq:14}, and the boundedness of $\big(\tilde{s}^{\eps_n},\tilde{z}^{\eps_n},\eps_n\tilde{\xi}^{\eps_n}\big)_n$, we directly compute that:
\begin{align}
\nonumber&-4C\eps_n^2(\eps_n w)^3(\tilde{s}^{\eps_n},\tilde{z}^{\eps_n},\tilde{\xi}^{\eps_n})
\left(w_{\xi^{i_0^n}}-w_{\xi^{j_0^n}}\right)(\tilde{s}^{\eps_n},\tilde{z}^{\eps_n},\tilde{\xi}^{\eps_n})\\
&\hspace{0.9em}+\eps_n^3v_z(\tilde{s}^{\eps_n},\tilde{z}^{\eps_n})\left[\lambda^{i_0^n,j_0^n}
-(1+\delta)(\partial_{i_0^n}\overline{w}-\partial_{j_0^n}\overline{w})(\tilde{s}^{\eps_n},\tilde{z}^{\eps_n},
\tilde{\rho}^{\eps_n})\right]+\circ(\eps_n^3)\geq 0.
\label{eq:16}
\end{align}
Using \reff{eq:15} and the non-negativity of $w$, this implies
 \b*
 0
 &\le&
 -4c_0\lambda^{i_0^n,j_0^n}\eps_n^2
  (\eps_nw)^3(\tilde{s}^{\eps_n},\tilde{z}^{\eps_n},\tilde{\xi}^{\eps_n})
 -\delta\lambda^{i_0^n,j_0^n}\eps_n^3
  v_z(\tilde{s}^{\eps_n},\tilde{z}^{\eps_n})
 +\circ(\eps_n^3)
 \\
 &\le& 
 -\delta\lambda^{i_0^n,j_0^n}\eps_n^3
  v_z(\tilde{s}^{\eps_n},\tilde{z}^{\eps_n})
  +\circ(\eps_n^3)\geq 0,
 \e*
which leads to a contradiction when $n$ goes to $+\infty$.

\vspace{0.35em}

{\bf Step $4$:} From \reff{eq:13} and Proposition \ref{prop}, we deduce that in our domain of interest, the function $\Psi^{\eps,\delta}$ is actually smooth, and is therefore a legitimate test function for the second order operator of the dynamic programming equation. We then obtain from the supersolution property of $v^\eps$ that
\begin{equation}
\left(\beta v^\eps-\mathcal L\Psi^{\eps,\delta}-\widetilde{U}(\Psi_x^{\eps,\delta})\right)(\tilde{s}^\eps,\tilde{x}^\eps,\tilde{y}^\eps)\geq 0.
\label{eq:17}
\end{equation}
Moreover, by Step $3$ and the continuity of $(s,z)\longmapsto\mathcal O(s,z)$ in Assumption \ref{assump.3}, the sequence $(\tilde{\xi}^\eps)_\eps$ is bounded. By classical results in the theory of viscosity solutions, there exists a sequence $\eps_n\rightarrow 0$ and some $\tilde\xi$ such that
 $$
 (s_n,z_n,\xi_n)
 :=
 (\tilde{s}^{\eps_n},\tilde{z}^{\eps_n},\tilde{\xi}^{\eps_n})
 \longrightarrow (s_0,z_0,\tilde\xi).
$$
Now recall that the function $w$ is smooth in this case and that the function $\Psi^{\eps,\delta}$ has exactly the form given in Section \ref{remainder}. By the remainder estimate in \reff{eq:17}, we obtain
\begin{equation*}
\frac12\eta(s_n,z_n)\abs{\sigma(s_n)\xi_n}^2+\frac12(1+\delta)
\Tr{\alpha\alpha^T(s_n,z_n)w_{\xi\xi}(s_n,z_n,\xi_n)}-\mathcal A\phi(s_n,z_n)+\mathcal R^\eps(s,z,\xi)\geq 0.
\end{equation*}
We still have no guarantee that $w$ is $C^2$ at $\widehat \xi$. Therefore, we carefully estimate the term involving $w_{\xi\xi}$. 
Indeed, the equation satisfied by $w$ yields
\begin{equation}
a(s_n,z_n)-\mathcal A\phi(s_n,z_n)+\delta\left(a(s_n,z_n)-\frac12\eta(s_n,z_n)\abs{\sigma(s_n)\xi_n}^2\right)+\mathcal R^\eps(s,z,\xi)\geq 0.
\label{eq:19}
\end{equation}
Notice that the estimate on the remainder of Section \ref{remainder} still hold true if terms involving $w_{\xi\xi}$ are replaced by means of the first corrector equation. Since the map $(s,z)\longmapsto a(s,z)$ is continuous, by Assumption \ref{assump.3}, and all derivatives of $\Phi^\eps$ vanish at the origin, we may send $\eps$ to $0$ in \reff{eq:19} and obtain
\begin{equation}
a(s_0,z_0)-\mathcal A\phi(s_0,z_0)+\delta\left(a(s_0,z_0)-\frac12\eta(s_0,z_0)|\sigma(s_0)\tilde\xi|^2\right)\geq 0.
\label{eq:20}
\end{equation}
Since $\tilde\xi$ is bounded uniformly in $\delta$, we let $\delta$ go to zero in \reff{eq:20} to obtain the desired result.
\qed

\subsection{Viscosity supersolution property}
In this section we will prove the following result.

\begin{Proposition}\label{super}
Under Assumptions \ref{assump2}, \ref{assump1} and \ref{assump.3}, $u_*$ is a viscosity supersolution of the second corrector equation.
\end{Proposition}

To prove this result, we start by two useful lemmas. For the first one, we recall that in the proof of the viscosity subsolution property in the previous section, we mentioned that the function $w$ is not $C^2$ in the whole space. We overcome this difficulty by using the fact that we only needed $w$ on a subset of $\mathbb R^d$ where it is actually smooth. However, in the proof of the viscosity supersolution property, we need $w$ to be defined on the whole space. Therefore, we mollify it and the following Lemma gives some useful properties satisfied by this mollified version of $w$. 

\vspace{0.35em}

Let $k:\mathbb R^d\rightarrow\mathbb R$ be a positive, even (i.e. $k(-x)=k(x)$ for all $x\in\mathbb R^d$), $C^\infty$ function with support in the closed unit ball of $\mathbb R^d$ and unit total mass. For all $m>0$, we define 
 \b*
 k^m(x)
 :=
 \frac{1}{m^d}k\Big(\frac{\xi}{m}\Big)
 &\mbox{and}&
 w^m(s,z,\xi)
 :=
 \int_{\mathbb R^d}k^m(\zeta)
                   w(\xi-\zeta)d\zeta
 -\int_{\mathbb R^d}k^m(\zeta)w(-\zeta)d\zeta.
 \e*

\begin{Lemma}\label{lemma.wh}
Let Assumption \ref{assump.3} hold. For any $m>0$, the function $w^m$ satisfies:

$\rm{(i)}$ $w^m$ is $C^2$, convex in $\xi$, and we have for all $0\leq i,j\leq d$
	$$
	w^m_{\xi_i}(s,z,\xi)
	=
	\int_{\mathbb R^d}k^m(\zeta)w_{\xi_i}(s,z,\xi-\zeta)d\zeta, 
	\text{ }
	w^m_{\xi_i\xi_j}(s,z,\xi)
	=\int_{\mathbb R^d}k^m(\zeta)w_{\xi_i\xi_j}(s,z,\xi-\zeta)d\zeta.
	$$
Moreover, $0\leq w^m(s,z,\xi)\leq Lv_z(s,z)\abs{\xi}.$

\vspace{0.35em}
$\rm{(ii)}$ $w^m$ is smooth in $(s,z)$, and satisfies the following estimates, uniformly in $m$, 
	\begin{align}\label{estim1h}
\nonumber&\left(\abs{w^m}+\abs{w^m_s}+\abs{w^m_{ss}}+\abs{w^m_z}+\abs{w^m_{sz}}+\abs{w^m_{zz}}\right)(s,z,\xi)\leq C(s,z)\left(1+\abs{\xi}\right)\\
\nonumber&\left(\abs{w^m_\xi}+\abs{w^m_{s\xi}}+\abs{w^m_{z\xi}}\right)(s,z,\xi)\leq C(s,z)\\
&\abs{w^m_{\xi\xi}(s,z,\xi)}\leq C(s,z)1_{\xi\in B(s,z)},
\end{align}
where $C(s,z)$ is a continuous function depending on the Merton value function and its derivatives, and $B(s,z)$ is small ball with continuous radius and center in $(s,z)$. 

\vspace{0.35em}
$\rm{(iii)}$ For every $0\leq i,j\leq d$ and every $(s,z,\xi)$
$$-\lambda^{i,j}v_z(s,z)+w^m_{\xi_i}(s,z,\xi)-w^m_{\xi_j}(s,z,\xi)\leq 0.$$

$\rm{(iv)}$ For every $(s,z,\xi)$, we have
 $$
 \frac12v_{zz}(s,z)\int_{\mathbb R^d}
                   k^m(x)\abs{\sigma(s)(\xi-\zeta)}^2d\zeta
 -\frac12\Tr{\alpha\alpha^T(s,z)w^m_{\xi\xi}(s,z,\xi)}+a(s,z)\leq 0.
 $$
\end{Lemma}

\proof
$\rm{(i)}$ The fact that $w^m$ is $C^2$ in $\xi$ is a classical result. Moreover, we have by definition $w^m(s,z,0)=0$ and by convexity of $w$, we have $w^m\geq w\geq 0$. The equalities for the derivatives of $w^m$ follow from the fact that $w$ is $C^1$ in $\xi$ and $C^2$ almost everywhere. Finally, $w^m$ inherits clearly the convexity and Lipschitz property of $w$.

\vspace{0.35em}
$\rm{(ii)}$ is clear by linearity of the convolution and Assumption \ref{assump.3}.
$\rm{(iii)}$ is again a consequence of the linearity of the convolution and the gradient constraints satisfied by $w$. Finally $\rm{(iv)}$ follows from the linearity of the convolution, the second corrector equation satisfied by $w$ and the formula for $w^m_{\xi\xi}$ given in $\rm{(i)}$. 
\qed

\vspace{0.4em}

We next constructs a useful function which plays a major role in our subsequent proof.

\begin{Lemma}\label{lemma.h}
For any $\delta\in(0,1)$, there exists $a^\delta>1$ and a function $h^\delta:\mathbb R^d\rightarrow [0,1]$ such that $h^\delta$ is $C^\infty$, $h^\delta=1$ on $\mathcal B_1(0)$ and $h^\delta=0$ on $B_{a^\delta}(0)^c$. Moreover, for any $1\leq i\leq d$ and for any $\xi\in\mathbb R^d$
 $$
 \abs{h^\delta_{\xi_i}(\xi)}
 \leq 
 \frac{\delta}{2L\overline{\lambda}}
 \1_{B_{a^\delta}(0)}(\xi),
 \text{ and }
 \abs{\xi}|h^\delta_{\xi\xi}|\leq C^*,
 $$
for some constant $C^*$ independent of $\delta$.
\end{Lemma}

\proof
Let $\phi$ be an even $C^\infty$ function on $\mathbb R_+$, whose support is in $(-1,1)$, such that $0\leq \phi\leq 1$ and $\int_{-1}^1\phi(t)dt=1$. For some $\alpha>0$ to be specified later, we define the following function $\tilde h^\alpha:\R_+\longrightarrow[0,1]$:
 $$
 \tilde h_\alpha(x)
 :=
 \int_{-1}^1H^\alpha(x-t)\phi(t)dt,
 \text{ where }
 H^\alpha(x)
 :=
 \1_{\{x\leq 2\}}
 +\Big(1-\alpha\ln\Big(\frac{x}{2}\Big)\Big)
  \1_{\{2<x\leq 2e^{1/\alpha}\}}.
 $$
Clearly, $\tilde h_\alpha$ is $C^\infty$, $\tilde h_\alpha=1$ on $[0,1]$ and to $\tilde h_\alpha=0$ for $[1+2e^{1/\alpha},\infty)$. Moreover, $\tilde h_\alpha$ is decreasing and its derivative clearly verifies for every $x\in\mathbb R$
\begin{equation}\label{h}
0\geq \tilde h_\alpha'(x)\geq -\frac\alpha2.
\end{equation}
We claim that
 \be\label{h1}
 \abs{x\tilde h_\alpha''(x)}\leq \frac{3\alpha}{4}
 &\mbox{for all}&
 x\in\mathbb R.
 \ee
Indeed, this inequality is obvious on $(-\infty,1]$, and we compute for $x\geq 1$, that
 $$
 0\leq x\tilde h_\alpha''(x)
 =
 \int_{-1}^{1}\frac{\alpha x}{(x-t)^2}
              \1_{\{2<x-t\leq 2e^{1/\alpha}\}}
              \phi(t)dt
 \le
 \int_{-1}^{1}\alpha\frac{3}{4}
              \phi(t)dt
 =
 \frac{3\alpha}{4}.
 $$
We now introduce the function:
 $$
 h^\delta(\xi)
 :=
 \tilde h_{\tilde\delta}\left(\abs{\xi}\right),
 ~~\tilde\delta:=\frac{\delta}{L\overline{\lambda}}
 ~~\xi\in\R^d.
 $$
Clearly $h^\delta$ is $C^\infty$, takes values in $[0,1]$, $h^\delta=1$ on $B_0(1)$, and $h^\delta=0$ on $B_{1+2e^{L/\delta}}(0)^c$. In particular, this provides the existence of $a^\delta \in[1,1+2e^{\frac L\delta}]$. Also
 $$
 h^\delta_{\xi_i}(\xi)
 =
 \tilde h_{\tilde\delta}'(\abs{\xi})\frac{\xi_i}{\abs{\xi}}.
 $$
Thus, by \reff{h} and the fact that $\tilde h_{\tilde\delta}$ and all its derivatives vanish on $(a^\delta,\infty)$:
 $$
 \abs{h^\delta_{\xi_i}(\xi)}
 \leq 
 \frac{\delta}{2L\overline\lambda}\1_{B_{a^\delta}(0)}(x).
 $$
Similarly, we have
 $$
 \abs{\xi}h^\delta_{\xi\xi}(\xi)
 =
 \Big(\abs{\xi}\tilde h_{\tilde\delta}''(\abs{\xi})
      -\frac12\tilde h_{\tilde\delta}'(\abs{\xi})
 \Big)
 \frac{\xi\xi^T}{\abs{\xi}^2}
 +\tilde h_{\tilde\delta}'(\abs{\xi})I_d.
 $$
Consequently, using \reff{h}, \reff{h1}, and the fact that $\delta\in(0,1)$, we have for some constant $C$ which only depends on the dimension $d$:
 $$
 \abs{\xi}\big|h^\delta_{\xi\xi}\big|
 \leq C\left(\frac{3\delta}{4L\overline\lambda}+\frac{\delta}{4L\overline\lambda}+\frac{\delta}{2L\overline\lambda}\right)\leq \frac{3C}{2L\overline\lambda}.
 $$
\qed

\vspace{0.4em}

\proof[Proposition \ref{super}]
Let $(s_0,z_0,\phi)\in(0,+\infty)^{d+1}\times C^2\left((0,+\infty)^{d+1}\right)$ be such that
\begin{equation}
(u_*-\phi)(s_0,z_0)<(u_*-\phi)(s,z),\text{ for all $(s,z)\in(0,+\infty)^d\times(0,+\infty)\backslash\left\{(s_0,z_0)\right\}$}.
\label{eq:200}
\end{equation}

By the definition of viscosity supersolutions, we need to show that 
\begin{equation}
\mathcal A\phi(s_0,z_0)-a(s_0,z_0)\geq 0.
\label{eq}
\end{equation}
We argue by contradiction and assume that 
\begin{equation}
\mathcal A\phi(s_0,z_0)-a(s_0,z_0)< 0.
\label{contra}
\end{equation}
Then by the continuity of $\phi$ and $a$, for some $r_0>0$, we will have
 \be\label{contra2}
 \mathcal A\phi(s,z)-a(s,z)\leq0
 &\mbox{on}&
 B_{r_0}(s_0,z_0)
 ~~\mbox{for some}~~r_0>0.
 \ee
{\bf Step $1$:} This first step is devoted to defining the test function we will consider in the sequel. First of all, we know from Lemma \ref{indep} that there exists a sequence $(s^\eps,z^\eps)$ which realizes the $\underline{\lim}$ for $\widehat u^\eps$, that is to say
$$(s^\eps,z^\eps)\underset{\eps\downarrow 0}{\longrightarrow} (s_0,z_0)\text{ and }\widehat u^\eps(s^\eps,z^\eps,0)\underset{\eps\downarrow 0}{\longrightarrow}u_*(s_0,z_0).$$

It follows then easily that $l^\eps_*:=\widehat u^\eps(s^\eps,z^\eps,0)-\phi(s^\eps,z^\eps)\underset{\eps\downarrow 0}{\longrightarrow}0,$ and
$$
(x^\eps,y^\eps)=\left(z^\eps-\ybf(s^\eps,z^\eps)\cdot \1_d,\ybf(s^\eps,z^\eps)\right)\underset{\eps\downarrow 0}{\longrightarrow}(x_0,y_0):=\left(z_0-\ybf(s_0,z_0)\cdot \1_d,\ybf(s_0,z_0)\right).
$$
We then choose $\eps_0$, depending on $z_0$, $s_0$ and $\phi$, such that for all $\eps\leq \eps_0$, we have
\begin{equation}
\abs{z^\eps-z_0}\leq \frac {r_0}{4},\ \abs{s^\eps-s_0}\leq\frac{r_0}{4},\ \abs{l^*_\eps}\leq 1.
\label{eq:1bis}
\end{equation}
By the continuity of $\phi$ and $v_z$, we may also introduce a constant $c_0>0$ such that 
 \be\label{eq:2bis}
 \underset{(s,z)\in\overline{ B}_{r_0/2}(s_0,z_0)}{\sup}
 \left\{\phi(s,z)+v_z(s,z)\right\}+3
 &\leq& 
 2c_0\left(\frac {r_0}{4}\right)^4.
 \ee
Notice that since $\phi$ and $v_z$ are continuous, the supremum on the compact set $\overline{\mathcal B}(s_0,z_0,r_0/2)$ above is indeed finite. This justifies the existence of $c_0$. We now define
 $$
 \varphi^\eps(s,z)
 :=
 \phi(s,z)-c_0\left(\abs{z-z^\eps}^4+\abs{s-s^\eps}^4\right).
 $$
Then, for all $\eps\leq \eps_0$ and for any $(s,z)\in\overline{B}_{r_0/2}(s_0,z_0)$, we have using \reff{eq:2bis} and \reff{eq:1bis}
\begin{align}\label{eq:3bis}
\nonumber
\varphi^\eps(s,z)+l^*_\eps+v_z(s,z)
=&\ 
\phi(s,z)+v_z(s,z)-c_0\left(\abs{z-z^\eps}^4+\abs{s-s^\eps}^4\right)+l^*_\eps
\\
\nonumber
\leq & \
2c_0\left(\frac {r_0}{4}\right)^4-3-c_0\left(\abs{z-z^\eps}^4+\abs{s-s^\eps}^4\right)+1
\\
\leq & \ c_0\left(2\left(\frac {r_0}{4}\right)^4-\abs{z-z^\eps}^4-\abs{s-s^\eps}^4\right)-2\leq -2,
\end{align}
whenever $(s,z)\in\partial B_{r_0/2}(s_0,z_0)$.
 
\vspace{0.35em}
Before defining our final test function, we provide another parameter. Let $\tilde\xi_0>0$ be greater than the diameter of the open bounded set $\mathcal O_0(s_0,z_0)$, and fix some $\xi^*\geq 1\vee\tilde\xi_0\vee\tilde\xi^*$ where
$$\tilde\xi^*:=\underset{(s,z)\in\overline{\mathcal B}(s_0,z_0,r_0/2)}{\sup}\left(\frac{1+2a(s,z)+\frac{C^*+\overline\lambda^{-1}}{2}\Tr{\alpha\alpha^T(s,z)}v_z(s,z)}{\frac{1}{2}\abs{\sigma}^2(-v_{zz})(s,z)}\right)^{1/2},$$
and $C^*$ is the constant appearing in Lemma \ref{lemma.h}(iv). Then, for any $\xi\in B_{\xi^*}(0)^c$, it follows from \reff{contra2} that for any $(s,z)\in\overline{B}_{r_0/2}(s_0,z_0)$:
\begin{equation}\label{contra3}
\frac{1}{2}\abs{\sigma\xi}^2(-v_{zz})(s,z)-\frac{C^*+\overline\lambda^{-1}}{2}\Tr{\alpha\alpha^T}v_z(s,z)-\mathcal A\phi(s,z)>1+2a(s,z)-\mathcal A\phi(s,z)>1.
\end{equation}

\vspace{0.35em}
Finally, for $\delta\in(0,1)$ and $m>0$, let $h^\delta$ be the function introduced in Lemma \ref{lemma.h}, $w^m$ the function introduced in Lemma \ref{lemma.wh}, and define the $C^2$ test function $\psi^{\eps,\delta,m}(s,x,y)$ and the corresponding $\hat\psi^{\eps,\delta,m}(s,z,\xi)$ by
 $$
 \hat\psi^{\eps,\delta,m}(s,z,\xi)
 :=
 v(s,z)-\eps^2\varphi^\eps(s,z)-\eps^2l^*_\eps
 -\eps^4(1-\delta)w^m(s,z,\xi)H(\xi),
 ~~H(\xi):=h^\delta\left(\frac{\xi}{\xi^*}\right).
 $$
{\bf Step $2$:}  In this step, we modify the test function once again
in order to recover the interior maximizer property. We first compute that:
 \be
 I^{\eps,\delta,m}(s,z,\xi)
 &:=&
 \eps^{-2}\big(\hat v^\eps-\hat\psi^{\eps,\delta,m}\big)(s,z,\xi)
 \nonumber\\
 &=&
 \widehat\phi^\eps(s,z)-\widehat u^\eps(s,z,\xi)+l^*_\eps-\eps^2\left[1-(1-\delta)H\left(\xi\right)\right]w^m(s,z,\xi).
\label{eq:6bis}
\ee
In particular, this implies
\begin{equation}
I^{\eps,\delta,m}(s^\eps,z^\eps,0)=\phi(s^\eps,z^\eps)-\widehat u^\eps(s^\eps,z^\eps,0)+l^*_\eps=0.
\label{eq:7bis}
\end{equation}
Furthermore, since $v^\eps\leq v$, we also have easily
\begin{equation}
I^{\eps,\delta,m}(s,z,\xi)\leq \varphi^\eps(s,z)+l^*_\eps+\eps^2(1-\delta)H(\xi)w^m(s,z,\xi).
\label{eq:8bis}
\end{equation}
Now, we use the fact that 
$$0\leq \delta\leq 1, \ v_z(s,z)>0, \ 0\leq H(\xi)\leq 1_{\abs{\xi}\leq a^\delta\xi^*} \text{ and }0\leq w^m(s,z,\xi)\leq Lv_z(s,z)\abs{\xi},$$ in \reff{eq:8bis} to obtain
 \be\label{eq:9bis}
 \nonumber 
 I^{\eps,\delta,m}(s,z,\xi)
 &\leq&
 \varphi^\eps(s,z)+l^*_\eps+\eps^2L(1-\delta)H(\xi)v_z(s,z)\abs{\xi}
 \\
 &\leq& \varphi^\eps(s,z)+l^*_\eps+\eps^2Lv_z(s,z)a^\delta\xi^*
 \;\leq\; 
 \varphi^\eps(s,z)+l^*_\eps+v_z(s,z),
 \ee
provided that $\eps\leq \eps^\delta:=(La^\delta\xi^*)^{-1/2}$.

\vspace{0.35em}
Define the set $Q_{s_0,z_0}:=\left\{(s,z,\xi),\ (s,z)\in\overline{B}_{r_0}(s_0,z_0)\right\}.$ Let us then distinguish two cases. First, we assume that $(s,z,\xi)\in\partial Q_{s_0,z_0}$ for every $\xi$. Then, if we take $\eps\leq \eps_0\wedge\eps^\delta$, using \reff{eq:3bis} and \reff{eq:9bis}, we obtain that for any $\xi$
\begin{equation}
I^{\eps,\delta,m}(z,\xi)\leq -2.
\label{eq:10bis}
\end{equation}
We assume next that $(s,z,\xi)\in \rm{int}\left(Q_{s_0,z_0}\right)$.  Then, once again for $\eps\leq \eps^\delta\wedge\eps_0$, we have
\begin{equation}
I^{\eps,\delta,m}(z,\xi)\leq \varphi^\eps(s,z)+l^*_\eps+v_z(s,z)\leq C(s_0,z_0)<+\infty,
\label{eq:11bis}
\end{equation}
for some constant $C(s_0,z_0)$ depending only on $\phi$, $s_0$ and $z_0$, since $(s,z)$ lies in a compact set, $l^*_\eps\leq 1$ and all the functions appearing on the right-hand side of \reff{eq:11bis} are continuous. This implies that 
$$\nu(\eps,\delta,m):=\underset{(s,z,\xi)\in Q_{s_0,z_0}}{\sup}\ I^{\eps,\delta,m}(s,z,\xi)<+\infty.$$

By definition, we can therefore for each $n\geq 1$ find $(\widehat s_n,\widehat z_n,\widehat \xi_n)\in\rm{int}\left(Q_{s_0,z_0}\right)$ such that
\begin{equation}
I^{\eps,\delta,m}(\widehat s_n,\widehat z_n,\widehat \xi_n)\geq \nu(\eps,\delta,m)-\frac{1}{2n}.
\label{eq:12bis}
\end{equation}
Since we have no guarantee that the maximizer above exists, we modify once more our test function. Let $f$ be an even smooth function such that $f(0)=1$, $f(x)=0$ if $\abs{x}\geq 1$ and $0\leq f\leq 1$. We then define the test function $\Psi^{\eps,\delta,m,n}(s,x,y)$ and the corresponding
 $$
 \hat\Psi^{\eps,\delta,m,n}(s,z,\xi)
 :=
 \hat\psi^{\eps,\delta,m}(s,z,\xi)
 -\frac{\eps^2}{n}f\big(|\xi-\widehat\xi_n|\big).
 $$
Consider then
 $$
 I^{\eps,\delta,m,n}(s,z,\xi)
 :=
 \eps^{-2}\left(\hat v^\eps-\hat\psi^{\eps,\delta,m,n}\right)(s,z,\xi)
 =
 I^{\eps,\delta,m}(s,z,\xi)+\frac1n f\left(|\xi-\widehat\xi_n|\right).
 $$
Notice now that for any $(s,z,\xi)\in Q_{s_0,z_0}$, we have using \reff{eq:12bis} 
\begin{equation}
I^{\eps,\delta,m,n}(\widehat s_n,\widehat z_n,\widehat\xi_n)=I^{\eps,\delta,m}(\widehat s_n,\widehat z_n,\widehat\xi_n)+\frac1n \geq \nu(\eps,\delta,m)+\frac{1}{2n}\geq I^{\eps,\delta,m}(s,z,\xi)+\frac{1}{2n}.
\label{eq:13bis}
\end{equation}
Moreover, by the definition of $f$, we have:
 \b*
 I^{\eps,\delta,m,n}(s,z,\xi)=I^{\eps,\delta,m}(s,z,\xi)
 &\mbox{if}&
 \xi\in B_1(\widehat \xi_n)^c
 ~~\mbox{and}~~
 (s,z,\xi)\in Q_{s_0,z_0}.
 \e*
This equality and \reff{eq:13bis} imply 
 \b*
 \underset{(s,z,\xi)\in Q_{s_0,z_0}}
 {\sup}\ I^{\eps,\delta,m,n}(s,z,\xi)
 &=&
 \underset{(s,z,\xi)\in \overline{B}_1(\widehat \xi_n)\cap Q_{s_0,z_0}}
 {\sup}\ I^{\eps,\delta,m,n}(s,z,\xi).
 \e*
Since $\overline{B}_1(\widehat \xi_n)\cap Q_{s_0,z_0}$ is a compact set, we deduce that there exists some $(s_n,z_n,\xi_n)\in Q_{s_0,z_0}$ which maximizes $I^{\eps,\delta,m,n}$. We now claim that we actually have $(s_n,z_n,\xi_n)\in\rm{int}(Q_{s_0,z_0})$. Indeed, we have by \reff{eq:7bis}
 $$
 I^{\eps,\delta,m,n}(s_n,z_n,\xi_n)
 \geq 
 I^{\eps,\delta,m,n}(s^\eps,z^\eps,0)
 \geq 
 I^{\eps,\delta,m}(s^\eps,z^\eps,0)
 =
 0,
 $$
and, by \reff{eq:10bis}, 
 \b*
 I^{\eps,\delta,m,n}(s,z,\xi)
 \leq 
 I^{\eps,\delta,m}(s,z,\xi)+\frac1n\leq -2+\frac1n<0,
 (s,z,\xi)\in\partial Q_{s_0,z_0}
 \mbox{ for }
 (s,z,\xi)\in\partial Q_{s_0,z_0}.
 \e*
By the viscosity subsolution property of $v^\eps$ at the point $(s_n,z_n,\xi_n)$, with corresponding $(s_n,x_n,y_n)$, it follows that
 \be\label{eq:14bis}
 \underset{0\leq i,j\leq d}{\min}
 \left\{\beta v^\eps-\mathcal  L\Psi^{\eps,\delta,m,n}
        -\widetilde U\left(\Psi_x^{\eps,\delta,m,n}\right)
        \ ,\    
        \Lambda^\eps_{i,j}
        \cdot
        (\Psi^{\eps,\delta,m,n}_x,\psi^{\eps,\delta,m,n}_y)
 \right\}
 &\le&
 0.
 \ee
{\bf Step $3$:} Our aim in this step is to show that for $\eps$ small enough and $n$ large enough, we have:
 \be\label{eq:15bis}
 D^{i,j}
 :=
 \Lambda^\eps_{i,j}
 \cdot
 \big(\Psi^{\eps,\delta,m,n}_x,\Psi^{\eps,\delta,m,n}_y\big)
 (s_n,x_n,y_n)
 \;>\;0
 &\mbox{for all}&
 0\leq i,j\leq d.
 \ee
We easily compute for $0\leq i\leq d$, with the convention that $y_0=x$ and $e_0=0$, that:
 \b*
 \Psi^{\eps,\delta,m,n}_{y_i}(s,x,y)
 &=&
 v_z(s,z)-\eps^2\varphi^\eps_z(s,z)
 -\eps^3(1-\delta)(w^mH)_\xi(s,z,\xi).(e_i-\ybf_z(s,z))
 \\
 &&
 -\eps^4(1-\delta)(w^m_zH)(s,z,\xi)-\frac\eps n\frac{f'(|\xi-\widehat\xi_n|)}{|\xi-\widehat\xi_n|}
 \xi\cdot(e_i-\ybf_z(s,z)).
 \e*
and
 \b*
 D^{i,j}
 &=&
 \eps^3\left[\lambda^{i,j}v_z(s_n,z_n)
             -(1-\delta)(w^mH)_{\xi}(s_n,z_n,\xi_n).(e_i-e_j)
       \right]
 -E^\eps-F^{\eps,n},
 \e*
where
\begin{align*}
&E^\eps:=\lambda^{i,j}\left[\eps^5(\phi_z(s_n,z_n)-4C(z_n-z^\eps)^3)+\eps^6(1-\delta)(w^mH)_{\xi}(s_n,z_n,\xi_n).(e_i-\ybf_z(s_n,z_n))\right]\\
&\hspace{2.6em}+\lambda^{i,j}\eps^7(1-\delta)(w^m_zH)(s_n,z_n,\xi_n),\\
&F^{\eps,n}:=\frac\eps n\frac{f'(|\xi_n-\widehat\xi_n|)}{|\xi_n-\widehat\xi_n|}\xi_n.\left(e_i-e_j+\lambda^{i,j}\eps^3(e_i-\ybf_z(s_n,z_n))\right).
\end{align*}
Recall that $\xi^*\geq 1$. Then, from Lemma \ref{lemma.h}, we have
 \b*
 0\leq H\leq 1,
 ~~\abs{H_\xi} \leq \frac{\sqrt{d}\delta}{2L}
 &\mbox{and}&
 H=0~~\mbox{for}~~\abs{\xi}\geq a^\delta\xi^*,
 \e*
we deduce
\begin{align}\label{eq:estimbis}
\nonumber\abs{E^\eps}\leq& \ \lambda^{i,j}\eps^5\left[\phi_z(s_n,z_n)+4c_0\abs{z_n-z^\eps}^3+\eps C(s_0,z_0)\left(\abs{w^m_\xi} H+w^m\abs{H_\xi}\right)(s_n,z_n,\xi_n)\right.\\
\nonumber&\left.+\eps^2Lv_z(s_n,z_n)\abs{\xi_n}1_{H(\xi_n)>0}\right]\\
\leq& \ C(s_0,z_0)\eps^5\left[1+\eps\left(1+a^\delta\right)+\eps^2a^\delta\xi^*\right],
\end{align}
for some constant $C(s_0,z_0)$ which can change value from line to line. Then we also have easily for some constant denoted $\mbox{Const}$, which can also change value from line to line
\begin{equation}
\abs{F^{\eps,n}}\leq \mbox{Const}\frac{\eps}{n}.
\label{eq:estim2}
\end{equation}

Let us now study the term
\begin{align*}
G^\eps:&=\lambda^{i,j}v_z-(1-\delta)(w^mH)_{\xi}.(e_i-e_j)\\
&=\lambda^{i,j}v_z-(1-\delta)(w^m_{\xi_i}-w^m_{\xi_j})H-(1-\delta)w^m(H_{\xi_i}-H_{\xi_j}),
\end{align*}
where we suppressed the dependence in $(s_n,z_n,\xi_n)$ for simplicity. Using Lemma \ref{lemma.wh}(iii) and the fact that $w^m$ and $H$ are positive, $w^m\leq Lv_z\abs{\xi}$ and $H=0$ for $\xi\geq a^\delta\xi^*$
\begin{align*}
G^\eps&\geq \lambda^{i,j}v_z-\lambda^{i,j}(1-\delta)v_z-(1-\delta)Lv_z\abs{\xi_n}\left(\abs{H_{\xi_i}}+\abs{H_{\xi_j}}\right)\\
&\geq\lambda^{i,j}v_z\left(\delta-\frac{La^\delta\xi^*}{\lambda^{i,j}}(1-\delta)\left(\abs{H_{\xi_i}}+\abs{H_{\xi_j}}\right)\right)\\
&\geq\lambda^{i,j}v_z\left(\delta-\frac{L\xi^*}{\overline\lambda}(1-\delta)\left(\abs{H_{\xi_i}}+\abs{H_{\xi_j}}\right)\right),
\end{align*}
since $a^\delta\geq 1$. But by (iii) of Lemma \ref{lemma.h}, we know that for all $0\leq i,j\leq d$, $\abs{H_{\xi_i}}\leq \frac{\delta}{2L\overline\lambda},$ from which we deduce
\begin{equation}
G^\eps\geq \lambda^{i,j}\delta^2v_z.
\label{eq:estim3}
\end{equation}

Finally, we have obtained
$$D^{i,j}\geq \lambda^{i,j}\delta^2v_z\eps^3-C(s_0,z_0)\eps^5\left(1+\eps\left(1+a^\delta\right)+\eps^2a^\delta\xi^*\right)
-\mbox{Const}\frac\eps n.$$

\vspace{0.35em}
Next, there is by hypothesis some constant $\widetilde C>0$ such that $v_z\geq \widetilde C$, and therefore there is a $\tilde\eps^\delta$ such that for all $\eps\leq \tilde\eps^\delta$:
$$C(s_0,z_0)\eps^5\left(1+\eps\left(1+a^\delta\right)+\eps^2a^\delta\xi^*\right)\leq \frac{\lambda^{i,j}\delta^2 \widetilde C\eps^3}{4}.$$
Then, for alll $n\ge N_{\eps,\delta}:=\frac{4\mbox{Const}}{\eps^2\widetilde C\delta^2}$, we have
$$
\mbox{Const}\frac\eps n\leq \frac{\lambda^{i,j}\delta^2 \widetilde C\eps^3}{4}.
$$
We then conclude that, for $\eps\leq \tilde\eps^\delta$ and $n\geq N_{\eps,\delta}$, we have
$D^{i,j}\geq\lambda^{i,j}\delta^2\widetilde C\eps^3/2>0$, and by the arbitrariness of $i,j=0,\ldots,d$, we deduce from \reff{eq:14bis} that
 \be
 J^{\eps,\delta,m,n}
 \;:=\;
 \frac{1}{\eps^2}\left[\beta v^\eps-\mathcal L\Psi^{\eps,\delta,m,n}
                       -\widetilde U\left(\Psi^{\eps,\delta,m,n}_x\right)
                 \right](s_n,x_n,y_n)\leq 0.
\label{eq:}
 \ee
{\bf Step $4$}: We now consider the remainder estimate. Using the general expansion result, we have
\begin{equation}
J^{\eps,\delta,m,n}=(-v_{zz})\frac{\abs{\sigma\xi_n}^2}{2}+\frac{1-\delta}{2}\Tr{\alpha\alpha^T(wH)_{\xi\xi}}-\mathcal A\phi+\mathcal R,
\label{eq:20bis}
\end{equation}

where
$$\mathcal R=\mathcal R_{\widetilde U}+\mathcal R_{\phi}+\mathcal R_{f}+\mathcal R_{wH},$$
and using the same calculations as in the remainder estimate of Section \ref{remainder}
\begin{align*}
&\abs{\mathcal R_{\widetilde U}}=\eps^{-2}\abs{\widetilde{U}(\Psi^{\eps,\delta,m,n}_x)-\widetilde U(v_z)-\eps^2\widetilde U^{'}(\phi_z)}\leq \mbox{Const}\left(\eps+\abs{\eps\xi_n}+\frac1n\right)\\
&\abs{\mathcal R_\phi}\leq \mbox{Const}\left[\abs{\eps\xi_n}+\abs{\eps\xi_n}^2\right]\\
&\abs{\mathcal R_f}\leq \frac{\mbox{Const}}{n}\\
&\abs{\mathcal R_{wH}}\leq \mbox{Const}\left(1+\abs{\eps\xi}^2\right)\left(\eps+\eps^2\abs{\xi_n}+\eps\abs{\xi_n}\right)1_{\abs{\xi_n}\leq a^\delta\xi^*}.
\end{align*}
We deduce from all these estimates that there is a $\widehat \eps^\delta$ and a $\tilde N_\eps$ such that if $\eps\leq \widehat \eps^\delta$ and $n\geq \tilde N_\eps$, we have
$$\abs{\mathcal R}\leq 1.$$

\vspace{0.35em}

{\bf Step $5$:} For $n$ greater than all the previously introduced $N$'s and $\epsilon$ smaller than all the ones previously introduced, we now show that $\abs{\xi_n}\leq \xi^*$.

\vspace{0.35em}
We argue by contradiction and we suppose that $\abs{\xi_n}>\xi^*$. Then we know that $w^m_{\xi\xi}(.,\xi_n)=0$. By \reff{eq:}, the expansion \reff{eq:20bis}, and the result of Lemma \ref{lemma.h}(iii),(iv), we see that:
\begin{align*}
(-v_{zz}(s_n,z_n))\frac{\abs{\sigma\xi_n}^2}{2}-\mathcal A\phi&\leq -\frac{1-\delta}{2}\Tr{\alpha\alpha^T(s_n,z_n)(w^mH)_{\xi\xi}(s_n,z_n,\xi_n)}-\mathcal R\\
&\leq -\frac{1-\delta}{2}\Tr{\alpha\alpha^T(s_n,z_n)\left(w^mH_{\xi\xi}+2w_{\xi}H_\xi^T\right)}+1\\
&\leq v_z(s_n,z_n)\frac{1-\delta}{2}\Tr{\alpha\alpha^T(s_n,z_n)}\left(\abs{\xi_n}\abs{H_{\xi\xi}}+\frac{\delta}{\overline\lambda}\right)+1\\
&\leq \frac{v_z(s_n,z_n)}{2}\Tr{\alpha\alpha^T(s_n,z_n)}\left(C^*+\frac{1}{\overline\lambda}\right)+1,
\end{align*} 
contradicting \reff{contra3}. Hence $(\xi_n)_n$ is bounded by $\xi^*$ (which does not depend on $\eps$, $\delta$, $m$ or $n$). In particular, this implies that the function $H$ applied to $\xi_n$ is always equal to $1$. Therefore, by the boundedness of $(s_n,z_n,\xi_n)$ and by classical results on the theory of viscosity solutions, there exists some $\bar\xi$ such that by letting $n$ go to $+\infty$ and then $\eps$ to $0$ (along some subsequence if necessary) in \reff{eq:20bis}, we obtain by using Lemma \ref{lemma.wh} (iv):
\begin{align}\label{eq:29}
\nonumber0&\geq -v_{zz}(s_0,z_0)\frac{|\sigma(s_0)\bar\xi|^2}{2}
+\frac{(1-\delta)}{2}\Tr{\alpha\alpha^T(s_0,z_0)w^m_{\xi\xi}(s_0,z_0,\bar\xi)}
-\mathcal A\phi(s_0,z_0)\\
\nonumber&\geq (1-\delta)a(s_0,z_0)-\delta v_{zz}(s_0,z_0)\frac{|\sigma(s_0)\bar\xi|^2}{2}-\mathcal A\phi(s_0,z_0)
\\
&\hspace{0.9em}+\frac{1-\delta}{2}v_{zz}(s_0,z_0)
\int_{\mathbb R^d} k^m(\zeta)\left(|\sigma(s_0)(\bar\xi-\zeta)|^2
                                    -|\sigma(s_0)\bar\xi|^2\right)d\zeta.
\end{align}
Now using the fact that the function $k^m$ is even, we have
$$
\int_{\mathbb R^d}k^m(\zeta)\left(|\sigma(s_0)(\bar\xi-\zeta)|^2
                            -|\sigma(s_0)\bar\xi|^2\right)d\zeta
=
\int_{\mathbb R^d}k^m(\zeta)\abs{\sigma(s_0)\zeta}^2d\zeta.
$$
Since $\int_{\mathbb R^d}k^m(\zeta)\abs{\sigma(s_0)\zeta}^2d\zeta\underset{m\rightarrow 0}{\longrightarrow} 0$, and $\bar\xi$ is uniformly bounded in $m$ and $\delta$, we can let $\delta$ and $m$ go to $0$ in \reff{eq:29} to obtain
$$\mathcal A\phi(s_0,z_0)-a(s_0,z_0) \geq 0,$$
which is the required contradiction to \reff{contra}, and completes the proof of the required result \reff{eq}.
\qed

\section{Wellposedness of the first corrector equation}
\label{sect:first corrector}

In this section, we collect the main proofs which allow us to obtain the wellposedness of the first corrector equation \reff{eq:corrector1}. Since the variables $(s,z)$ are frozen in this equation, we simplify the notations by suppressing the dependence on them.

\vspace{0.3em}
Recall that since the set $C$ is bounded, convex and closed, the supremum in the definition of the convex function $\delta_C$ is always attained at the boundary $\partial C$. Moreover $0\in \text{$\rm{int}$}(C)$, we may find two constants $L,L'>0$ such that
\begin{equation}\label{norm}
L'\abs{\rho}\leq \delta_C(\rho)\leq L\abs{\rho}.
\end{equation}

\subsection{Uniqueness and comparison}
\label{sect:comparison w}

\no {\bf Proof of Theorem \ref{prop.uni}}
Fix some $(\eps,\nu)\in(0,1)\times(0,+\infty)$ and 
define for $(\rho,y)\in\mathbb R^d\times\mathbb R^d$
$$
\overline{w}^\eps(\rho,y)
:=
(1-\eps) w_1(\rho)-w_2(y),\ \phi_\nu(\rho,y)=\frac{1}{2\nu}\abs{\rho-y}^2.
$$
Since $w_1$ is a viscosity subsolution of \reff{eq:corrector1}, 
then its gradient takes values in $C$ in the viscosity sense, which implies that $w_1$ is $L-$Lipschitz. Then, for $y\neq 0$:
\begin{align*}
\overline{w}^\eps(\rho,y)-\phi_\nu(\rho,y)
&=(1-\eps)(w_1(\rho)-w_1(y))
-\frac{1}{2\nu}\abs{\rho-y}^2+(1-\eps) w_1(y)-w_2(y)
\\
&\leq (1-\eps) L\abs{\rho-y}-\frac{1}{2\nu}\abs{\rho-y}^2
 +(1-\eps) w_1(y)-w_2(y)
\\
&\leq \frac{(1-\eps)^2L^2\nu}{2}+(1-\eps) w_1(y)-w_2(y)
\\
&=\frac{(1-\eps)^2L^2\nu}{2}
  +\delta_C(y)\left((1-\eps)\frac{w_1(y)}{\delta_C(y)}
                     -\frac{w_2(y)}{\delta_C(y)}
              \right).
\end{align*}
By the growth conditions on $w_1,w_2$, together with \reff{norm}, this implies that:
$$
\underset{\abs{(\rho,y)}\rightarrow+\infty}{\lim}
\overline{w}^\eps(\rho,y)-\phi_\nu(\rho,y)=-\infty.
$$
Then, the difference $\overline{w}^\eps-\phi_\nu$ has a global maximizer $(\rho^{\eps,\nu},y^{\eps,\nu})\in\mathbb R^d\times\mathbb R^d$ satisfying the lower bound
\begin{equation}
(\overline{w}^\eps-\phi_\nu)(\rho^{\eps,\nu},y^{\eps,\nu})
\geq (\overline{w}^\eps-\phi_\nu)(0,0)=0.
\label{eq:finite}
\end{equation}
By the Crandall-Ishii Lemma (see Theorem $3.2$ in \cite{cil}), it follows that for any $\eta>0$, there exist symmetric positive matrices $X$ and $Y$ such that
\begin{align}
\label{ishii}
\nonumber &\left(D_\rho\phi_\nu(\rho^{\eps,\nu},y^{\eps,\nu}),X\right)
=\left(\frac{\rho^{\eps,\nu}-y^{\eps,\nu}}{\nu},X\right)
\in\overline{J}^{2,+}((1-\eps)w_1)(\rho^{\eps,\nu})\\
&\left(-D_y\phi_\nu(\rho^{\eps,\nu},y^{\eps,\nu}),Y\right)
=\left(\frac{\rho^{\eps,\nu}-y^{\eps,\nu}}{\nu},Y\right)
\in\overline{J}^{2,-}w_2(y^{\eps,\nu}),
\end{align}
and 
\b*
\left(
  \begin{array}{ c c }
     X & 0 \\
     0 & -Y
  \end{array} \right)\leq A+\eta A^2,
&\mbox{with}&
A:=D^2\phi_\nu(\rho^{\eps,\nu},y^{\eps,\nu})=\frac1\nu\left(
\begin{array}{ c c }
I_d & -I_d \\
-I_d & I_d
\end{array} \right).
\e*
The above matrix inequality 
directly implies that $X\leq Y$.
We now use \reff{ishii} to arrive at
$$
\frac{\rho^{\eps,\nu}-y^{\eps,\nu}}{(1-\eps)\nu}
\in\overline{J}^{1,+}w_1(\rho^{\eps,\nu}).
$$
In addition, since $w_1$ is a viscosity subsolution, $Dw_1\in C$ in the viscosity sense. This implies that for all $0\leq i,j\leq d$,
$$
-\lambda^{j,i}\leq \frac{\rho^{\eps,\nu}_i
-y^{\eps,\nu}_i}{(1-\eps)\nu}-\frac{\rho^{\eps,\nu}_j
-y^{\eps,\nu}_j}{(1-\eps)\nu}\leq \lambda^{i,j}.
$$
Since $\eps\in(0,1)$, we have
\b*
-\lambda^{j,i}< \frac{\rho^{\eps,\nu}_i
-y^{\eps,\nu}_i}{\nu}-\frac{\rho^{\eps,\nu}_j
-y^{\eps,\nu}_j}{\nu}< \lambda^{i,j},
&\mbox{and therefore}&
\frac{\rho^{\eps,\nu}-y^{\eps,\nu}}{\nu}
\in\text{$\rm{int}$}(C).
\e*
Consequently, it follows from \reff{ishii} and the viscosity subsolution and supersolution of $w_1$ and $w_2$ that
\b*
-\frac{\abs{\sigma \rho^{\eps,\nu}}^2}{2}
-\frac{1}{2(1-\eps)}\Tr{\bar{\alpha}\bar{\alpha}^T X}
+a_1
&\le\; 0\;\le&
-\frac{\abs{\sigma y^{\eps,\nu}}^2}{2}
-\frac12\Tr{\bar{\alpha}\bar{\alpha}^T Y}+a_2.
\e*
Since $X\le Y$, this provides:
\begin{equation}
(1-\eps)a_1-a_2
\leq 
\frac12\Tr{\bar{\alpha}\bar{\alpha}^T(X-Y)}
+(1-\eps)\frac{\abs{\sigma \rho^{\eps,\nu}}^2}{2}
-\frac{\abs{\sigma y^{\eps,\nu}}^2}{2}
\leq \frac{\abs{\sigma \rho^{\eps,\nu}}^2}{2}
-\frac{\abs{\sigma y^{\eps,\nu}}^2}{2}.
\label{eq:ab}
\end{equation}
We now show that $\big(y^{\eps,\nu},\rho^{\eps,\nu}\big)_\nu$ remains bounded as $\nu$ tends to zero.
\begin{itemize} 
\item We argue by contradiction, assuming to the contrary that
\b*
\abs{y^{\eps,\nu_n}} \longrightarrow \infty
&\mbox{for some sequence}&
\nu_n\longrightarrow\infty.
\e*
Since $w_1$ is Lipschitz, this implies that
\begin{align*}
\big(\overline{w}^\eps-\phi_{\nu_n}\big)
(\rho^{\eps,\nu_n},y^{\eps,\nu_n})
\leq\frac{((1-\eps)L)^2\nu_n}{2}
+\delta_C(y^{\eps,\nu_n})\frac{1-\eps}{\delta_C(y^{\eps,\nu_n})}\big(w_1(y^{\eps,\nu_n})-w_2(y^{\eps,\nu_n})\big).
\end{align*}
Arguing as in the beginning of this proof, we see that
$(\overline{w}^\eps-\phi_{\nu_n})(\rho^{\eps,\nu_n},y^{\eps,\nu_n})\longrightarrow-\infty$, contradicting \reff{eq:finite}.
\item Similarly, using the normalization $w_1(0)=0$, we have
\begin{align*}
\big(\overline{w}^\eps-\phi_{\nu_n}\big)(\rho^{\eps,\nu_n},y^{\eps,\nu_n})
\leq(1-\eps) L\abs{\rho^{\eps,\nu_n}}-w_2(y^{\eps,\nu_n})
-\frac{\abs{\rho^{\eps,\nu_n}-y^{\eps,\nu_n}}^2}{\nu_n}.
\end{align*}
Since $(y^{\eps,\nu_n})_n$ was just shown to be bounded, we see that 
$\abs{\rho^{\eps,\nu_n}}\longrightarrow\infty$ implies $(\overline{w}^\eps-
-\phi_{\nu_n})(\rho^{\eps,\nu_n},y^{\eps,\nu_n})\longrightarrow-\infty$,
contradiction. 
\end{itemize}
By standard techniques from the theory of viscosity solutions, we may then construct a $\rho^\eps\in\mathbb R^d$ and a sequence  $(\tilde\nu_n)_{n\geq 0}$ converging to zero such that
$
\big(\rho^{\eps,\tilde\nu_n},y^{\eps,\tilde\nu_n}\big)
\longrightarrow(\rho^\eps,\rho^\eps)$, as $n\to\infty$. Passing to the limit in \reff{eq:ab} along this sequence, we see that $(1-\eps)a_1-a_2\leq 0,$ which implies that $a_1\le a_2$ by the arbitrariness of $\eps \in (0,1)$.
\qed

\vspace{0.4em}

The following uniqueness result is an immediate consequence.

\begin{Corollary}
There is at most one $a\in\mathbb R$ such that \reff{eq:corrector1} 
has a viscosity solution $\overline{w}$ satisfying the growth condition $\overline{w}(\rho)/\delta_C(\rho)\longrightarrow 1$, as $|\rho|\to\infty$.
\end{Corollary}

\begin{Remark}
{\rm{ (i)
In the context of \cite{hyn}, $C$ is just the closed unit ball, then it is clear that $\delta_C(\rho)=\abs{\rho}$, and the growth condition in the previous result reduces to $\overline{w}(\rho)/|\rho|\longrightarrow 1$, as $|\rho|\to\infty$.

\vspace{0.3em}
(ii) In the one-dimensional case $d=1$, we directly compute that $\delta_C(\rho)=\lambda^{1,0}\rho^++\lambda^{0,1}\rho^-$.
So the above growth condition is the sharpest one for the 
explicit solution of \reff{eq:corrector1} given in Section $4$ of \cite{st}.}}
\qed
\end{Remark}

\subsection{Optimal control approximation of the first corrector equation}\label{sect:existence w}

In ergodic control, it is standard to introduce an approximation by a sequence of infinite horizon standard control problems with a vanishing
discount factor $\eta>0$:
\begin{eqnarray}
\nonumber
&&\underset{0\leq i,j\leq d}{\max}\max\left\{-\frac{\abs{\sigma \rho}^2}{2}
-\frac12\Tr{\bar{\alpha}\bar{\alpha}^T D^2\overline{w}^\eta(\rho)}
+\eta \overline{w}^\eta(\rho),\right.\\
&&\hspace{150pt}
\left.
-\lambda^{i,j}
+ \frac{\partial \overline{w}^\eta}{\partial \rho_i}(\rho)
- \frac{\partial \overline{w}^\eta}{\partial \rho_j}(\rho) \right\}=0,
\label{eq:corrector22}
\end{eqnarray}
together with the growth condition
 \be\label{eq:growth}
\underset{\abs{\rho}\rightarrow+\infty}{\lim}\frac{\overline{w}^\eta(\rho)}{\delta_C(\rho)}=1.
 \ee

\vspace{0.4em}
We first state a comparison result for the equation \reff{eq:corrector22}-\reff{eq:growth}. The proof is omitted as it is very similar to that of Theorem \ref{prop.uni}

\begin{Proposition}
Let $w_1,w_2$ be respectively a viscosity subsolution and a viscosity supersolution of \reff{eq:corrector22}. Assume further that
$$
\underset{\abs{\rho}\rightarrow \infty}{\overline\lim}\ \frac{w_1(\rho)}{\delta_C(\rho)} 
\leq1\leq \underset{\abs{\rho}\rightarrow \infty}{\underline\lim}\ \frac{w_1(\rho)}{\delta_C(\rho)}.
$$
Then, $w_1\leq w_2$. In particular, there is at most one viscosity solution of \reff{eq:corrector22}. 
\end{Proposition}

The next result states the existence of a unique solution of the approximating control problem.

\begin{Proposition}\label{exist}
For every $\eta\in(0,1]$, there is a unique viscosity solution $\overline{w}^\eta$ of \reff{eq:corrector22}-\reff{eq:growth}. Moreover, $\overline{w}^\eta$ is $L$-Lipschitz {\rm{(}}with a constant $L$ independent of $\eta${\rm{)}}, and we have the following estimate
\begin{equation}
\left(\delta_C(\rho)-K_1\right)^+\leq \overline{w}^\eta(\rho)\leq\frac{K_2}{\eta}
+\delta_C(\rho),\ \rho\in\mathbb R^d.
\label{eq:estim}
\end{equation}
\end{Proposition}

\proof In view of the previous comparison result, we establish existence of a viscosity solution by an application of Perron's method, which requires to find appropriate sub and supersolutions. The remaining properties are immediate consequences. We then introduce:
\b*
\underline\varpi(\rho)
:=
\left(\delta_C(\rho)-K_1\right)^+,
\mbox{ and }
\overline\varpi(\rho)
:=
\frac{K_2}{\eta}
+\1_{\{\delta_C(\rho)<1\}}
\frac{\delta_C(\rho)^2}{2}
+\1_{\{\delta_C(\rho)\geq 1\}}\left(\delta_C(\rho)-\frac12\right).
\e*
$\rm{(i)}$ We first prove that we may choose $K_1$ so that $\underline\varpi$ is a viscosity subsolution of the equation \reff{eq:corrector22} satisfying the growth condition \reff{eq:growth}. Since $\underline\varpi$ has linear growth, $\delta_C$ is Lipschitz and vanishes at $0$, we may choose $K_1>0$ such that 
$\underline \varpi(\rho)\leq \abs{\sigma \rho}^2/2$ for all $\rho\in\mathbb R^d$.

\vspace{0.3em}
Moreover, $\underline\varpi$ is convex and has a gradient in the weak sense, which takes values in $C$, by definition of the support function. Then, for all $\rho_0\in\mathbb R^d$, and $(p,X)\in J^{2,+}\underline\varpi(\rho_0)$, we have $X\geq 0$, $p\in C$, and therefore:
\begin{align*}
&\underset{0\leq i,j\leq d}{\max}\max\left\{-\frac{\abs{\sigma \rho_0}^2}{2}
-\frac12\Tr{\bar{\alpha}\bar{\alpha}^T X}+\eta \underline\varpi(\rho_0),-\lambda^{i,j}+ p_i- p_j \right\}\\
&\hspace{100pt}\leq \underset{0\leq i,j\leq d}{\max}\max\left\{-\frac{\abs{\sigma \rho_0}^2}{2}
+\eta \underline\varpi(\rho_0),-\lambda^{i,j}+ p_i- p_j \right\}
\;\leq\; 0.
\end{align*}
Hence, $\underline\varpi$ is a viscosity subsolution.

\vspace{0.3em}
$\rm{(ii)}$ We next prove that $\overline\varpi$ is a viscosity supersolution of the equation \reff{eq:corrector22} satisfying the growth condition \reff{eq:growth}, for a convenient choice of $K_2$. Consider arbitrary $\rho_0\in\mathbb R^d$ and $(p,X)\in J^{2,-}\overline\varpi(\rho_0)$. 

\vspace{0.3em}
\underline{Case 1}: $\delta_C(\rho_0)<1$, then $\rho_0$ must be bounded. Moreover, by definition, $\overline\varpi$ is Lipschitz. Hence,  
$p\in C$ and $X\geq 0$. In particular, $p$ is bounded and by the
definition of $J^{2,-}\overline{\varpi}(\rho_0)$ so is $X$. 
We then have
$$
-\frac{\abs{\sigma \rho_0}^2}{2}-\frac12\Tr{\bar{\alpha}\bar{\alpha}^T X}
+\eta \overline\varpi(\rho_0)\geq -\frac{\abs{\sigma \rho_0}^2}{2}-\frac12\Tr{\bar{\alpha}\bar{\alpha}^T X} +K_2.
$$
Since $\rho_0$ and $X$ are bounded,
we may also choose $K_2$ large enough so that the above is positive.
This implies the supersolution property in that case.

\vspace{0.3em}
\underline{Case 2}: $\delta_C(\rho_0)\geq 1$, then $\overline\varpi$ has a weak gradient which, by definition of the support function, takes values in $\partial C$, i.e. at least one of the gradient constraints in \reff{eq:corrector22} is binding. This implies that the supersolution property is satisfied.
\qed

\vspace{0.4em}

We next establish that $\overline{w}^\eta$ is convex, by following the PDE argument of \cite{hyn}. Notice that this property would have been easier to prove if the probabilistic representation of Remark \ref{rem.rep} was known to be valid.

\begin{Lemma}\label{prop.conv}
$\overline{w}^\eta$ is convex and therefore is twice differentiable Lebesgue almost everywhere.
\end{Lemma}

\proof 1.
Let $\eps\in(0,1)$, $\rho^0,\rho^1\in\R^d$, $\bar\rho:=(\rho^0+\rho^1)/2$, and let us first prove that
 \b*
 \ell^{\eps}(\rho^0,\rho^1)
 :=
 (1-\eps) \overline{w}^\eta(\bar\rho)
 - \big(\overline{w}^\eta(\rho^0)+\overline{w}^\eta(\rho^1)\big)/2
 \longrightarrow-\infty
 &\mbox{as}&
 |(\rho^0,\rho^1)|\to\infty.
 \e*
Let $(\rho^0_n,\rho^1_n)\in\mathbb R^d\times\mathbb R^d$ be such that $\abs{\rho^0_n}+\abs{\rho^1_n}\longrightarrow\infty.$ Denote $\bar\rho_n:=(\rho^0_n+\rho^1_n)/2$. For large $n$, we have $\delta_C(\rho^0_n)+\delta_C(\rho^1_n)>0$, and using the convexity of $\delta_C$ we also have
 \b*
 \Delta_n:=
 \frac{\ell^{\eps}(\rho^0_n,\rho^1_n)}
      {\delta_C(\rho^0_n)+\delta_C(\rho^1_n)}
 &=&
 (1-\eps)\frac{\overline{w}^\eta(\bar\rho_n)}
               {\delta_C(\rho^0_n)+\delta_C(\rho^1_n)}
 -\frac12\sum_{i=0,1}
   \frac{\delta_C(\rho^i_n)}
               {\delta_C(\rho^i_n)+\delta_C(\rho^{1-i}_n)}
   \frac{\overline{w}^\eta(\rho^i_n)}{\delta_C(\rho^i_n)}
 \\
 &\le&
 \frac{1-\eps}{2}\frac{\overline{w}^\eta(\bar\rho_n)}
                      {\delta_C(\bar\rho_n)}
 -\frac12\sum_{i=0,1}
   \frac{\delta_C(\rho^i_n)}
               {\delta_C(\rho^i_n)+\delta_C(\rho^{1-i}_n)}
   \frac{\overline{w}^\eta(\rho^i_n)}{\delta_C(\rho^i_n)}.
 \e*
Consider first the case where $\big(\delta_C(\bar\rho_n)\big)_n$ is bounded, which is equivalent to the boundedness of $\big(\abs{\bar\rho_n}\big)_n$ by \reff{norm}. Then it is clear by the growth property \reff{eq:growth} that $\Delta_n\longrightarrow-\frac12<0.$
Similarly, if $\delta_C(\bar\rho_n)\longrightarrow\infty$, we see that
$\limsup_{n\to\infty}\Delta_n\le [(1-\eps)-1]/2<0$. In both case this proves the required result of this step.

\vspace{0.3em}
2. For $\theta=(\rho^0,\rho^1,y^0,y^1)\in\R^{4d}$, set $\bar\rho:=(\rho^0+\rho^1)/2$, $\bar y:=(y^0+y^1)/2$, and define:
 $$
 \varpi^{\eps}(\theta)
 :=
 (1-\eps)\overline{w}^\eta(\bar\rho)
 -\big(\overline{w}^\eta(y^0)+\overline{w}^\eta(y^1)\big)/2,
 ~~
 \phi_n(\theta)
 :=
 n\left(\abs{\rho^0-y^0}^2+\abs{\rho^1-y^1}^2\right)/2.
 $$
Since $\overline{w}^\eta$ is $L$-Lipschitz, we have
 \b*
 (\varpi^{\eps}-\phi_n)(\theta)
 &=&
 (1-\eps)\big[(\overline{w}^\eta(\bar\rho)
               -\overline{w}^\eta(\bar y)
         \big]
 -\phi_n(\theta)
 +\ell^\eps(y^0,y^1)
 \;\le\;
 \;\le\; 
 \ell^{\eps}(y^0,y^1)+L^2/n.
 \e*
By the first step, this shows that $(\varpi^{\eps}-\phi_n)(\theta)\longrightarrow-\infty$ as $|\theta|\longrightarrow\infty$, and that there is a global maximizer $\theta_n:=(\rho^0_n,\rho^1_n,y^0_n,y^1_n)$ of the difference $\varpi^{\eps}-\phi_n$.
Using then Crandall-Ishii's lemma and arguing exactly as in the proof of Proposition \ref{prop.uni} (see also the proof of Lemma $3.7$ in \cite{hyn}), we obtain that for all $(\rho^0,\rho^1)$
 \be
 (1-\eps) \overline{w}^\eta(\bar\rho)
 -\big(\overline{w}^\eta(\rho^0)
      +\overline{w}^\eta(\rho^1)
  \big)/2
 &\le&
 (1-\eps) \overline{w}^\eta(\bar\rho_n)
 -\big(\overline{w}^\eta(y^0_n)
       +\overline{w}^\eta(y^1_n)
  \big)/2
 \nonumber\\
 &\le&  \big(2\abs{\sigma\bar\rho_n}^2
             -|\sigma y^0_n|^2-|\sigma y^1_n|^2
        \big)/(4\eta).
 \label{eq:crandall}
 \ee
Following the same arguments as in the proof of Proposition \ref{prop.uni}, we next show that the sequence $(\theta_n)_n$ is bounded, and that there is a subsequence which converges to some $\theta_\eps:=(\rho^0_\eps,\rho^0_\eps,y^0_\eps,y^1_\eps)$. The averages $\bar\rho_\eps$ and $\bar y_\eps$ are introduced similarly. Passing to the limit along this subsequence in \reff{eq:crandall}, we obtain for all $(\rho^0,\rho^1)$:
 \b*
 \eta\left((1-\eps) \overline{w}^\eta(\bar\rho)
           -\big[\overline{w}^\eta(\rho^0)+\overline{w}^\eta(\rho^1)\big]/2
     \right)
 &\le&  
 \left(2\abs{\sigma\bar\rho_\eps}^2
       -\abs{\sigma \rho^0_{\eps}}^2
       -\abs{\sigma \rho^1_{\eps}}^2
 \right)
 \;\le\;
 0.
 \e*
The proof is completed by sending $\eps$ to $0$.
\qed

\vspace{0.4em}

As a consequence of the convexity of $\overline{w}^\eta$, we have the following result.

\begin{Lemma}\label{lemma.p}
There is a constant $M>0$ independent of $\eta$ such that $J^{1,-}\overline{w}^\eta(\rho)\subset\partial C$, for all
$\rho\in B_0(M)^c$.
\end{Lemma}

\proof
Let $K_2$ be the constant in Proposition \ref{exist}.
Since  $\delta_C$ has linear growth,
we can choose  $M$ large enough so that
$$
K_2+\delta_C(y)<\frac{\abs{\sigma y}^2}{2}\ \ 
\text{ for all}\ \ \abs{y}\geq M.
$$
Fix $\rho\in\mathbb R^d$ satisfying $\abs{\rho}\geq M$. 
Then, by the convexity of $\overline{w}^\eta$,  
$p\in J^{1,-}\overline{w}^\eta(\rho)$ if and only if $p$ belongs to the 
subdifferential of $\overline{w}^\eta$ at $\rho$. 
This implies, in particular, that $p\in C$. Furthermore, 
we have $(p,0)\in J^{2,-}\overline{w}^\eta(\rho)$.  Then,
the supersolution property of $\overline{w}^\eta$ yields
$$
\underset{0\leq i,j\leq d}{\max}\max\left\{-\frac{\abs{\sigma \rho}^2}{2}
+\eta \overline{w}^\eta(\rho),-\lambda^{i,j}+ p_i- p_j \right\}\geq 0.
$$
By the definition of $M$, we have $
-\frac{\abs{\sigma \rho}^2}{2}
+\eta \overline{w}^\eta(\rho)
\leq -\frac{\abs{\sigma \rho}^2}{2}+ K+\delta_C(\rho)<0
$, and therefore $\underset{0\leq i,j\leq d}{\max}\left\{-\lambda^{i,j}+ p_i- p_j \right\}\geq 0$. 
Since $p\in C$, this means that this quantity is actually equal to zero, 
implying that $p\in\partial C$.
\qed

\vspace{0.4em}

The following result is similar to Lemma $3.10$ in \cite{hyn}.

\begin{Lemma}\label{hess}
For Lebesgue almost every $\rho\in\mathbb R^d$, we have
$0\leq D^2\overline{w}^\eta(\rho)\leq\frac1\eta\No{\sigma\sigma^T}.$
\end{Lemma}

\proof
We fix $\eps\in(0,1)$ and $z\in\mathbb R^d$ such that 
$0<\abs{z}<M$  which we will send to zero later. Set
$$
\psi(\rho):=(1-\eps)(\overline{w}^\eta(\rho+z)
+\overline{w}^\eta(\rho-z))-2\overline{w}^\eta(\rho).
$$
Since $\overline{w}^\eta$ is  Lipschitz,
$$
\psi(\rho)\leq 2(1-\eps) L\abs{z}-2\overline{w}^\eta(\rho)
\leq 2(1-\eps) LM-2\overline{w}^\eta(\rho).
$$
In view of  the growth condition \reff{eq:growth}, as $\rho$
approaches to infinity
$\psi(\rho)$ tends to minus infinity.
Let
$$
\beta(\rho_1,\rho_2,\rho_3):=(1-\eps)(\overline{w}^\eta(\rho_1+z)
+\overline{w}^\eta(\rho_2-z))-2\overline{w}^\eta(\rho_3),
$$
and
$$
\phi_\nu(\rho_1,\rho_2,\rho_3):=\frac{1}{2\nu}\left(\abs{\rho_1-\rho_3}^2
+\abs{\rho_2-\rho_3}^2\right).
$$
We then have, again by the Lipschitz property of $\overline{w}^\eta$,
that
\begin{eqnarray*}
(\beta-\phi_\nu)(\rho_1,\rho_2,\rho_3)&=&
(1-\eps)\left(\overline{w}^\eta(\rho_1+z)
-\overline{w}^\eta(\rho_3+z)+\overline{w}^\eta(\rho_2+z)
-\overline{w}^\eta(\rho_3-z)\right)
\\
&&
+\psi(\rho_3)-\frac{1}{2\nu}\left(\abs{\rho_1-\rho_3}^2+\abs{\rho_2-\rho_3}^2\right)\\
&\leq& L^2\nu +\psi(\rho_3) \ \to\ -\infty, \quad {\mbox{as}}\ \abs{\rho} \to \infty.
\end{eqnarray*}

Hence there exists $(\rho_1^\nu,\rho_2^\nu,\rho_3^\nu)$ which maximizes $\beta-\phi_\nu$. By applying Crandall-Ishii's lemma, for every $\rho>0$, we can find symmetric matrices $(X,Y)\in\mathbb S^{2d}\times\mathbb S^{d}$ such that
\begin{align*}
\left(D_{\rho_1}\phi_\nu(\rho_1^\nu,\rho_2^\nu,\rho_3^\nu),
D_{\rho_2}\phi_\nu(\rho_1^\nu,\rho_2^\nu,\rho_3^\nu),X\right)
& \in\overline{J}^{2,+}(1-\eps)\left(\overline{w}^\eta(\rho_1^\nu+z)+\overline{w}^\eta(\rho_2^\nu-z)\right)\\
\left(-D_{\rho_3}\phi_\nu(\rho_1^\nu,\rho_2^\nu,\rho_3^\nu),Y\right)
&\in\overline{J}^{2,-}2\overline{w}^\eta(\rho_3^\nu),
\end{align*}
and
\begin{equation}
\label{eq:matrix2}
\left(
  \begin{array}{ c c c}
     X & 0 &0\\
     0 & -Y & 0\\
     0 & 0 & 0
  \end{array} \right)\leq A+\rho A^2,
  \end{equation}
where
$$A:=\left(
  \begin{array}{ c c }
     D^2\phi_\nu(\rho^{\nu}_1,\rho^{\nu}_2,\rho_3^\nu) & 0 \\
     0 & 0 
  \end{array} \right)=\frac1\nu\left(
  \begin{array}{ c c c c}
     I_d & 0& -I_d & 0\\
     0 & I_d & -I_d & 0\\
     -I_d & -I_d & 2I_d & 0\\
     0 & 0 & 0 & 0
  \end{array} \right).$$

For $(X_1,X_2,X_3)\in\mathbb S^d\times
\mathbb R^{d\times d}\times\mathbb S^d$,
set  $X:=\left(\begin{array}{ c c }
     X_1 & X_2 \\
     X_2^T & X_3
  \end{array} \right)$. Then, \reff{eq:matrix2} implies that
  $$
  X\leq\left(\begin{array}{ c c }
     -Y & 0 \\
     0 & 0
  \end{array} \right),
  $$
  and in particular, $\Tr{X_1+X_3-Y}\leq 0$.
  
\vspace{0.35em}  
We  directly calculate that
\begin{align*}
&\left(\frac{\rho_1^\nu-\rho_3^\nu}{\nu},X_1\right)
\in\overline{J}^{2,+}\left(\rho_1\mapsto (1-\eps)\left( \overline{w}^\eta(\rho_1+z)
+\overline{w}^\eta(\rho_2^\nu-z)\right)\right)_{\rho_1=\rho_1^\nu}\\
&\left(\frac{\rho_2^\nu-\rho_3^\nu}{\nu},X_3\right)
\in\overline{J}^{2,+}\left(\rho_2\mapsto (1-\eps)\left( \overline{w}^\eta(\rho_1^\nu+z)
+\overline{w}^\eta(\rho_2-z)\right)\right)_{\rho_2=\rho_2^\nu}\\
&\left(\frac{\rho_1^\nu+\rho_2^\nu-2\rho_3^\nu}{\nu},Y\right)
\in\overline{J}^{2,-}2\overline{w}^\eta(\rho_3^\nu).
\end{align*}

Since $\overline{w}^\eta$ is a viscosity subsolution of \reff{eq:corrector2}, 
we deduce that for all $0\leq i,j\leq d$
$$-\lambda^{j,i}\leq \frac{\rho_{1,i}^\nu-\rho_{3,i}^\nu}{(1-\eps)\nu}
 -\frac{\rho_{1,j}^\nu-\rho_{3,j}^\nu}{(1-\eps)\nu}\leq \lambda^{i,j},$$
and   
$$-\lambda^{j,i}\leq \frac{\rho_{2,i}^\nu-\rho_{3,i}^\nu}{(1-\eps)\nu} 
-\frac{\rho_{2,j}^\nu-\rho_{3,j}^\nu}{(1-\eps)\nu}\leq \lambda^{i,j}.$$

From this we deduce that for all $0\leq i,j\leq d$
$$-\lambda^{j,i}\leq \frac{\rho_{1,i}^\nu+\rho_{2,i}^\nu-2\rho_{3,i}^\nu}{2\eps\nu} -\frac{\rho_{1,j}^\nu
+\rho_{2,j}^\nu-2\rho_{3,j}^\nu}{2(1-\eps)\nu}\leq \lambda^{i,j}.$$

Since $\eps\in(0,1)$, this implies that $\frac{\rho_{1}^\nu
+\rho_{2}^\nu-2\rho_{3}^\nu}{2\nu}\in\text{$\rm{int}$}(C)$. 
Also $\frac{\rho_{1}^\nu+\rho_{2}^\nu-2\rho_{3}^\nu}{2\nu}
\in \overline{J}^{2,-}\overline{w}^\eta(\rho_3^\nu)$. 
Given that $\overline{w}^\eta$ is a viscosity supersolution, we deduce that
$$-\frac{\abs{\sigma \rho_3^\nu}^2}{2}-\frac12\Tr{\bar{\alpha}\bar{\alpha}^T\frac Y2}
+\eta \overline{w}^\eta(\rho_3^\nu)\geq 0.$$

Since $\overline{w}^\eta$ is  a viscosity subsolution, 
$$
-\frac{\abs{\sigma (\rho_1^\nu+z)}^2}{2}
-\frac{1}{2(1-\eps)}\Tr{\bar{\alpha}\bar{\alpha}^TX_1}
+\eta \overline{w}^\eta(\rho_1^\nu+z)\leq 0,
$$
and
$$
-\frac{\abs{\sigma (\rho_2^\nu-z)}^2}{2}
-\frac{1}{2(1-\eps)}\Tr{\bar{\alpha}\bar{\alpha}^TX_3}
+\eta \overline{w}^\eta(\rho_2^\nu-z)\leq 0.
$$
Summing up the last three inequalities, we obtain that for all $\rho\in\mathbb R^d$,
\begin{align}\label{ineq}
\nonumber
&(1-\eps)\left(\overline{w}^\eta(\rho+z)+\overline{w}^\eta(\rho-z)\right)
-2\overline{w}^\eta(\rho)\leq(1-\eps)\left(\overline{w}^\eta(\rho_1^\nu+z)
+\overline{w}^\eta(\rho_2^\nu-z)\right)-2\overline{w}^\eta(\rho_3^\nu)\\
\nonumber
&\hspace{50pt}
\leq \frac{1}{2\eta}\Tr{\bar{\alpha}\bar{\alpha}^T(X_1+X_3-Y)}
+\frac{(1-\eps)}{2\eta}\left(\abs{\sigma (\rho_1^\nu+z)}^2
+\abs{\sigma (\rho_2^\nu-z)}^2\right)-\frac{\abs{\sigma \rho_3^\nu}^2}{\eta}\\
&\hspace{50pt}
\leq \frac{1}{2\eta}\left(\abs{\sigma (\rho_1^\nu+z)}^2
+\abs{\sigma (\rho_2^\nu-z)}^2\right)-\frac{\abs{\sigma \rho_3^\nu}^2}{\eta}.
\end{align}

We argue as in the proofs of Proposition \ref{prop.uni} and Lemma \ref{prop.conv},
to show that $(\rho_1^\nu,\rho_2^\nu,\rho_3^\nu)$ is bounded and that there is a subsequence 
and a vector $\rho^*$ such that they all converge to $\rho^*$ along this subsequence. 
Using this result in \reff{ineq}, we obtain
\begin{align*}
\eta\left((1-\eps)\left(\overline{w}^\eta(\rho+z)+\overline{w}^\eta(\rho-z)\right)
-2\overline{w}^\eta(\rho)\right)&\leq \frac12\left(\abs{\sigma (\rho^*+z)}^2
+\abs{\sigma (\rho^*-z)}^2\right)-\abs{\sigma \rho^*}^2\\
&\leq \No{\sigma\sigma^T}\abs{z}^2,
\end{align*}
where we used the mean value theorem and 
the fact that the Hessian matrix of the map
$\rho \rightarrow {\abs{\sigma \rho}^2}/{2}$ is $\sigma\sigma^T$.

\vspace{0.35em}
We finally let $\eps$ go to $0^+$, divide the inequality by 
$\abs{z}^2$ and let $\abs{z}$ go to zero to obtain the result, 
since we know that the second derivative of $\overline{w}^\eta$ 
exists almost everywhere.
\qed

\vspace{0.4em}

\begin{Corollary}\label{cor}
\begin{itemize}
	\item[$\rm{(i)}$] $\overline{w}^\eta \in C^{1,1}(\mathbb R^d)$.
	\item[$\rm{(ii)}$] The set
	$$
	\mathcal O_\eta:=\left\{\rho\in\mathbb R^d,\ D\overline{w}^\eta(\rho)\in\text{$\rm{int}$}(C)\right\}.
	$$
	 is open and bounded independently of $0<\eta\leq 1$.
	\item[$\rm{(iii)}$] $\overline{w}^\eta\in C^\infty(\mathcal O_\eta)$.
	\item[$\rm{(iv)}$] There exists a constant $L>0$ independent of $\eta$ such that 
	$$
	0\leq D^2\overline{w}^\eta(\rho)\leq L,\ \rho\in\mathcal O_\eta.
	$$ 
\end{itemize}
\end{Corollary}

\begin{Remark}\label{rem.rep2}
{\rm{
Now that we have obtained some regularity for the function $\overline w^\eta$, we could hope to prove that it is indeed the value function of an infinite horizon stochastic control problem by means of a verification argument. However, the problem here is that the verification argument needs to identify an optimal control, which we can not do in the present setting. Indeed, we do not know whether the solution to the reflexion problem on the free boundary defined by PDE \reff{eq:corrector2} has a solution, since we know nothing about the regularity of the free boundary. Notice that such a regularity was assumed in \cite{mrt} in order to construct an optimal control.}}
\end{Remark}

\proof
The first result is a simple consequence of the previous Lemma \ref{hess}. 
For the second one, the fact that $\mathcal O_\eta$ is bounded independently of 
$\eta$ follows from Lemma \ref{lemma.p}, and it is open because it is the inverse image 
of the open interval $(-\infty,0)$ by the continuous application 
$$
\rho\mapsto \underset{0\leq i,j\leq d}{\max}\left\{-\lambda^{i,j}
+\frac{\partial \overline{w}^\eta}{\partial \rho_i}(\rho)-\frac{\partial \overline{w}^\eta}{\partial \rho_j}(\rho)\right\}.
$$

Then the third result follows from classical regularity results for
 linear elliptic PDEs (see Theorem $6.17$ in \cite{gil}), 
 since we have on $\mathcal O^\eta$
$$
-\frac{\abs{\sigma \rho}^2}{2}
-\frac12\Tr{\bar{\alpha}\bar{\alpha}^TD^2\overline{w}^\eta(\rho)}
+\eta \overline{w}^\eta(\rho)=0.
$$

Finally, by convexity of $\overline{w}^\eta$, 
we have for any $\xi$ such that $\abs{\xi}=1$ 
and for any $\rho\in\mathcal O_\eta$
\begin{align*}
\xi^TD^2\overline{w}^\eta(\rho)\xi
&\leq C_0\Tr{\bar{\alpha}\bar{\alpha}^TD^2\overline{w}^\eta(\rho)}
=2C_0\left(\eta \overline{w}^\eta(\rho)-\frac{\abs{\sigma \rho}^2}{2}\right)
\leq 2C_0(K_2+\delta_C(\rho)), 
\end{align*}
and the result is a consequence of the fact that $\mathcal O_\eta$ is bounded.
\qed

\vspace{0.4em}

In the following result, we extend Proposition $3.12$ of \cite{hyn} to our context and show that $\overline{w}^\eta$ is 
characterized by its values in $\overline{\mathcal O}_\eta$. 
The proof is very similar to the one given
in \cite{hyn}. We provide 
it here for the convenience of the reader.

\begin{Proposition}\label{prop.ext}
We have for all $\rho\in\mathbb R^d$, $
\overline{w}^\eta(\rho)
=\underset{y\in\overline{\mathcal O}_\eta}{\inf}\left\{\overline{w}^\eta(y)
+\delta_C(\rho-y)\right\}.
$
\end{Proposition} 

\proof
First of all, since $\overline{w}^\eta\in C^{1,1}(\mathbb R^d)$,
it is easy to see by the mean value Theorem and the fact that 
$D\overline{w}^\eta\in C$ that for all $(\rho,y)\in\mathbb R^d\times\mathbb R^d$
$$
-\delta_C(y-\rho)\leq \overline{w}^\eta(\rho)-\overline{w}^\eta(y)\leq\delta_C(\rho-y).
$$

Hence we have $\overline{w}^\eta(y)+\delta_C(\rho-y)\geq \overline{w}^\eta(\rho).$ This implies that for $\rho\in\overline{\mathcal O}_\eta$ we have $$\underset{y\in\overline{\mathcal O}_\eta}{\inf}\left\{\overline{w}^\eta(y)+\delta_C(\rho-y)\right\}=\overline{w}^\eta(\rho).$$

It remains to show that the result also holds for $\rho\in\overline{\mathcal O}_\eta^c$. Notice first that by convexity of $\overline{w}^\eta$, the minimum in the formula is necessarily achieved on $\partial \mathcal O_\eta$. Define then
$$\widetilde{w}^\eta(\rho):=\underset{y\in\partial\mathcal O_\eta}{\inf}\left\{\overline{w}^\eta(y)+\delta_C(\rho-y)\right\}.$$

Let us then consider the following PDE
\begin{align*}
\nonumber&\underset{0\leq i,j\leq d}{\max}\left\{-\lambda^{i,j}+\frac{\partial v}{\partial \rho_i}(\rho)-\frac{\partial v}{\partial \rho_j}(\rho)\right\}=0, \ \rho\in\overline{\mathcal O}_\eta^c\\
&v(\rho)=\overline{w}^\eta(\rho),\ \rho\in\partial\mathcal O_\eta.
\end{align*}

It is easy to adapt the proof of Proposition \ref{prop.uni} to obtain that this PDE satisfies a comparison principle. Since $\overline{w}^\eta$ clearly solves it, it is the unique solution. Then by the usual properties of inf-convolutions, it is clear that the function $\widetilde{w}^\eta$ which is convex, has a gradient in the weak sense which is in $C$. This implies that $\widetilde{w}^\eta$ is a subsolution of the above PDE. Moreover, since $\delta_C$ also has a gradient in the weak sense which is in $\partial C$ (by convexity), using the fact that by properties of the inf-convolution of convex functions, the subgradient of $\widetilde{w}^\eta$ at any point is contained in a subgradient of $\delta_C$, we have that $\widetilde{w}^\eta$ dominates every subsolution of the PDE which coincides with $\overline{w}^\eta$ on $\partial \mathcal O_\eta$. By usual results of the theory of viscosity solutions (see \cite{cil}) this proves that $\widetilde{w}^\eta$ is also a supersolution. By uniqueness, we obtain the desired result.
\qed

\vspace{0.4em}

\begin{Corollary}\label{cor.estim}
For Lebesgue almost every $\rho\in\overline{\mathcal O}_\eta^c$, $
D^2\overline{w}^\eta(\rho)=0.
$ Hence, the second derivative of $\overline{w}^\eta$ is 
bounded almost everywhere, independently 
of $\eta$ in the whole space $\mathbb R^d$.
\end{Corollary}

\proof
Let $\rho\in\overline{\mathcal O}_\eta^c$. 
By Proposition \ref{prop.ext}, there exists $y\in\partial\mathcal O_\eta$ 
such that $\overline{w}^\eta(\rho)=\overline{w}^\eta(y)+\delta_C(\rho-y)$. 
We then have for any $z\in\mathbb R^d$
\begin{align*}
\frac{\overline{w}^\eta(\rho+z)+\overline{w}^\eta(\rho-z)
-2\overline{w}^\eta(\rho)}{\abs{z}^2}&=\frac{\overline{w}^\eta(\rho+z)
+\overline{w}^\eta(\rho-z)-2\overline{w}^\eta(y)-2\delta_C(\rho-y)}{\abs{z}^2}\\
&\leq \frac{\delta_C(\rho+z-y)
+\delta_C(\rho-z-y)-2\delta_C(\rho-y)}{\abs{z}^2}.
\end{align*} 

Now, the function $\rho\rightarrow\delta_C(\rho)$ 
is continuous and clearly $C^\infty$ almost everywhere (actually except on $d(d+1)$ 
hyperplanes which therefore have Lebesgue measure $0$). 
More than that, this function is piecewise linear, which implies that its 
Hessian is null almost everywhere. By letting $\abs{z}$ go to $0$ above, 
we obtain the desired result. The last result is now a simple consequence of Corollary \ref{cor}$\rm{(iv)}$.
\qed

\vspace{0.4em}
Then, the uniform estimates obtained above allow us to prove our main existence result in Theorem \ref{prop}.

\vspace{0.4em}

\no {\bf Proof of Theorem \ref{prop}}
By the uniform estimates of Corollary \ref{cor.estim}, we may follow the arguments of Section $4.2$ in \cite{hyn} to construct a strictly positive sequence $(\eta_n)_{n\geq 0}$ converging to zero, 
$\overline a\in\mathbb R$ and $\overline w\in C^{1,1}(\mathbb R^d)$ 
such that
$$
\underset{n\rightarrow +\infty}{\lim}\eta_n\overline{w}^{\eta_n}(\rho^{\eta_n})
=\overline a,\text{ and }\overline{w}^{\eta_n}\underset{n\rightarrow+\infty}{\longrightarrow}\overline{w}
\text{ in $C^1_{loc}(\mathbb R^d)$},
$$
where $\rho^{\eta_n}$ is a global minimizer of $\overline{w}^{\eta_n}$.
Corollary \ref{cor} implies that the limiting function $\overline{w}$ satisfies the required properties.
\qed

\end{document}